\documentclass[english]{article}
\usepackage[T1]{fontenc}
\usepackage[latin9]{inputenc}
\usepackage{geometry}
\usepackage{verbatim}
\usepackage{amsmath}
\usepackage{amssymb}
\usepackage{esint}

\makeatletter
\newcommand{\lyxaddress}[1]{
\par {\raggedright #1
\vspace{1.4em}
\noindent\par}
}

\makeatother

\makeatother

\usepackage{babel}

\begin{document}

\title{Indefinite integrals of Bessel and Struve functions with half-integer
indices, and incomplete gamma functions with integer indices, all
multiplied by $Exp(-ax^{2})$ and divided by powers}

\author{Jack C. Straton}

\maketitle

\lyxaddress{Department of Physics, Portland State University, Portland, OR, 97207-0751,
straton@pdx.edu}
\begin{abstract}
Indefinite integrals are found for 63 half-integer Bessel and Struve
functions, and incomplete gamma functions with integer indices, each
multiplied by $Exp(-ax^{2})$ and divided by powers. A series solution
is given for the individual terms (of any inverse power) in such functions,
split into even and odd portions. Eight integrals are given involving
these series. Three infinite sets of auxiliary integrals are also
given.
\end{abstract}
\vspace{2pc}
 \textit{Keywords}: Macdonald functions; modified spherical Bessel
functions of the second kind; reduced Bessel functions; spherical
Bessel functions; modified spherical Bessel functions of the first
kind; spherical Bessel functions of the first kind; Neumann Functions;
spherical Bessel functions of the second kind; Meijer G-functions;
Struve functions; and Incomplete Gamma functions\\
 \\
 MSC classes: 44A20; 44A30; 81Q99; 33C10; 33C20; 33C60;33B20; 34B27;
30E20; 28-02

\section{Introduction}

Bessel functions play a prominent role as solutions to problems in
mathematics, physical sciences, and engineering. Researchers use them
for a broad range of electromagnetic problems at both the classical
and quantum levels, and viscoelastic and dielectric relaxation processes
are modeled using Bessel and Struve functions. Likewise, incomplete
gamma functions appear in a vast array of problems in mathemastics,
science, and engineering. Therefore integrals of such functions provide
the launching point for innumerable branches of inquiry. While tables
of \emph{definite} integrals of these functions are extensiove, there
very few \emph{indefinite} integrals containing any sort of Bessel
function in either Prudnikov, Brychkov, and Marichev~\cite{PBM1}
(Section 1.12.2) or Gradshteyn and Ryzhik~\cite{GR5} (Section 5.5),
and fewer combine these with exponentials and/or powers. This paper
addresses that derth.

Specifically, one would like to find indefinite integrals of half-integer
Macdonald functions (also called the {}``modified spherical Bessel
function of the second kind'' or {}``reduced Bessel~functions''),
multiplied by $e^{-ax^{2}}$ and divided by powers,

\begin{equation}
\int e^{-ax^{2}}x^{\frac{1}{2}-h}K_{n+\frac{1}{2}}\left(bx\right)\, ds\:.\label{eq:exp(s^2) K(s)}\end{equation}
 One might hope to expand the Macdonald function and integrate this
term by term, but one finds no indefinite integrals of the form

\begin{equation}
\int\frac{e^{-a^{2}x^{2}-bx}}{x^{h}}\, dx\:,\label{eq:exp(-x^2-x)/x^h}\end{equation}

\noindent in either Prudnikov, Brychkov, and Marichev \cite{PBM1}
(Section 1.3.3) or Gradshteyn and Ryzhik \cite{GR5} (Section 2.32),
even for $h=1$ or $2$, no less for a general inverse power. Several
alternative failed attempts at integrating (\ref{eq:exp(s^2) K(s)})
are given in Appendix B of a prior paper \cite{stra26a}. 

It turns out that some \emph{sums} of inverse powers, each having
the above form with different values of \emph{h}, can be integrated,
particularly those consisting of half-integer Bessel and Struve functions
and incomplete gamma functions with integer indices. This paper generates
and proves a large subset of these indefinite integrals.

\subsection{One context}

This investigation grew out of my development of a one-range addition
theorem that has no infinite second series~\cite{stra26a}, which
I used to give a series solution to a transition amplitude that Cheshire
set up in 1964 \cite{Cheshire}, but could not fully solve (and which
had defied reduction to analytic form ever since): the Fourier transform
of a product of Slater orbitals

\begin{equation}
\int d^{3}x_{1}x_{1}^{-1}e^{-\eta_{1}x_{1}}x_{12}^{-1}e^{-\eta_{2}x_{12}}e^{-i\mathbf{k}\cdot\mathbf{x}_{1}}\:,
\label{eq:kyy-1}\end{equation}

\noindent 
in which $\mathbf{x}_{12}=\mathbf{x}_{1}-\mathbf{x}_{2}$ is short-hand
form for shifted coordinates, $x_{2}$ is a coordinate variable external
to the integration, and $\mathbf{k}$ is a momentum variable within
the plane wave paired with the integration variable. Application of
the addition theorem gives a series of integrals

\begin{equation}
\int_{\eta_{2}}^{\eta_{1}}dx\, b^{n+\frac{1}{2}}x^{-2j-2m+3n+\frac{1}{2}}exp\left(\frac{ix^{2}k\cdot b}{\eta_{2}^{2}-\eta_{1}^{2}}\right)K_{n+\frac{1}{2}}\left(bx\right)\:\label{eq:s-integral}\end{equation}

\noindent in which $0\leq j\leq n$, and likewise for \emph{m}. The
terms with $g\equiv-2j-2m+2n\geq0$ can be integrated in terms of
the error~function via~\cite{PBM1} (p. 139 No. 1.3.2.5 for even
non-negative powers and p. 140 No. 1.3.2.6 for odd positive powers
(noting that the exponentials have different coefficients: $e^{-a^{2}x^{2}}$
versus $e^{-ax^{2}}$, respectively)). However, powers more negative
than $g$ within an expansion of the Macdonald function cannot so
be integrated, nor do they seem to be tabled anywhere.

This derth is somewhat surprising since many of these integrals can
be found.

\subsection{Notation}

The notation for these integrals in Sections 2-5 is 

\begin{equation}
\begin{array}{ll}
V\left(h,\, n,\, a,\, b,\, x\right)=\sqrt{\frac{2}{\pi}}\int e^{-ax^{2}}b^{n+\frac{1}{2}}x^{h+\frac{1}{2}}\frac{1}{2}G_{0,2}^{2,0}\left(\frac{x^{2}b^{2}}{4}|\begin{array}{c}
\frac{1}{2}\left(n+\frac{1}{2}\right),\frac{1}{2}\left(-n-\frac{1}{2}\right)\end{array}\right)\, dx\\
=\sqrt{\frac{2}{\pi}}\int e^{-ax^{2}}b^{n+\frac{1}{2}}x^{h+\frac{1}{2}}K_{n+\frac{1}{2}}(xb)\, dx\end{array}\:\:.\label{eq:V(h,n,x)}\end{equation}
where $n$ corresponds to the index of the Macdonald function $K_{n+\frac{1}{2}}(xb)$
(also appearing in the power of the prepending coefficient $\sqrt{\frac{2}{\pi}}b^{n+\frac{1}{2}}$
that makes the analytical results somewhat cleaner) and $h$ corresponds
to the prepending power $x^{h+\frac{1}{2}}$ throughout all that follows.
A Meijier G-function form of the Macdonald function is also given
in each of the following on the first line, and the explicit expansion
in inverse powers will be given on the third line of each expression below.

Not every integral of this general type can be found either through
conventional methods presented in Sections 2-5, or through symbolic
computation that one may use to check one's results, as discussed
in Section 6. So Section 7 introduces a series solution for the generic
integral (\ref{eq:exp(-x^2-x)/x^h}) whose notation $\Upsilon\left(h,\, a,\, b,\, x\right)$
has no mention of \emph{n}. The series results are made cleaner
by avoiding the prepending coefficient $\sqrt{\frac{2}{\pi}}b^{n+\frac{1}{2}}$,
that made the analytical-result integrals, above, cleaner.

Section 8 gives indefinite integrals of incomplete gamma functions
having integer indices, each multiplied by $Exp(-ax^{2})$ and divided
by powers and Sections 9-12 give such indefinite integrals of other sorts
of half-integer Bessel and Struve functions.

\section{Integrals with $h_{max}=2$ }

Consider the integral (which is written in three equivalent forms),

\textbf{Theorem 1}

\begin{equation}
\begin{array}{ll}
V\left(-1,\,1,\, a,\, b,\, x\right)=\sqrt{\frac{2}{\pi}}\int e^{-ax^{2}}b^{n+\frac{1}{2}}x^{-1+\frac{1}{2}}\frac{1}{2}G_{0,2}^{2,0}\left(\frac{b^{2}x^{2}}{4}|\begin{array}{c}
\frac{3}{4},-\frac{3}{4}\end{array}\right)\, dx\\
=\sqrt{\frac{2}{\pi}}\int e^{-ax^{2}}b^{1+\frac{1}{2}}x^{-1+\frac{1}{2}}K_{\frac{3}{2}}(xb)\, dx\\
=b\int\frac{1}{x}e^{-xb-ax^{2}}\left(\frac{1}{xb}+1\right)\, dx=-\frac{e^{-ax^{2}-bx}}{x}-\sqrt{a}\sqrt{\pi}e^{\frac{b^{2}}{4a}}\text{erf}\left(\frac{2ax+b}{2\sqrt{a}}\right)\end{array}\:\:\label{eq:V(-1,1,x) Theorem 1}\end{equation}
\textbf{Proof:}

One may prove this theorem using integration by parts with \begin{equation}
u=e^{-ax^{2}-bx}\label{eq:u}\end{equation}
so that \begin{equation}
\text{du}=e^{-ax^{2}-bx}(-2ax-b)\:.\label{eq:du}\end{equation}
Writing the first term in parentheses in (\ref{eq:V(-1,1,x) Theorem 1})
with a generic power as

\begin{equation}
\text{dv}=\frac{1}{x^{\upsilon}}\:,\label{eq:dv}\end{equation}
then \begin{equation}
v=\frac{x^{1-\upsilon}}{1-\upsilon}\quad\left[\upsilon\neq1\right]\label{eq:v}\end{equation}
and \begin{equation}
\int u\,\text{dv}=uv-\int v\, du=\frac{x^{1-\upsilon}e^{-ax^{2}-bx}}{1-\upsilon}-\int e^{-ax^{2}-bx}\left(-\frac{2ax^{2-\upsilon}}{1-\upsilon}-\frac{bx^{1-\upsilon}}{1-\upsilon}\right)\, dx\:\label{eq:uv-intvdu}\end{equation}

In the present case one has $\upsilon=2$ so that \cite{PBM1} (p.
140 No. 1.3.3.17)

\begin{equation}
\begin{array}{ll}
V\left(-1,\,1,\, a,\, b,\, x\right)=\int e^{-bx-ax^{2}}\left(\frac{1}{x^{2}}+\frac{b}{x}\right)\, dx=-\frac{e^{-ax^{2}-\text{bx}}}{x}-\int2ae^{-ax^{2}-bx}\, dx-\int(\frac{b}{x}-\frac{b}{x})e^{-ax^{2}-bx}\, dx\\
=-\frac{e^{-ax^{2}-bx}}{x}-\sqrt{a}\sqrt{\pi}e^{\frac{b^{2}}{4a}}\text{erf}\left(\frac{2ax+b}{2\sqrt{a}}\right)\end{array}\:,\label{eq:V(-1,1,x) proof}\end{equation}
thus completing the proof. $\square$

Note that nothing in this proof restricts the parameters \emph{a}
and \emph{b} in any way. Indeed, the contextual integral (\ref{eq:s-integral})
has an imaginary value for \emph{a} that posed no numerical or analytical
difficulties in that prior work \cite{stra26a}. The same lack of
restriction applies throughout the following proofs, except where
explicilty noted, as where we note that integrals with any non-negative
integer power of \emph{x} multiplying $e^{-ax^{2}-bx}$ in Prudnikov,
Brychkov, and Marichev, \cite{PBM1}(p. 140 No. 1.3.3.16,17) require
\emph{a} to be real (and positive) but there is no such restriction
on \emph{b}. Indeed, in Section 9 we have \emph{b} take on negative
values, and imaginary values in Sections 10-12.

This proof likewise does not restrict the variable \emph{x}, though
one usually assumes the integration variable \emph{x} to be a real
variable. (That is, $\int f\left(x\right)\, dx=\int_{0}^{x}f\left(x'\right)\, dx'$.)
It is in the \emph{application} of these formulas that care must be
taken. For instance, when parameters \emph{a} and \emph{b} are real,
while $x\rightarrow0$, the first term of the result, above, goes
to $-\infty$ while the second goes to zero if \emph{b} is also zero
\cite{Erf(0)}. If $x\rightarrow\infty$, the first term goes to $0$ while the
second goes to $-\sqrt{a}\sqrt{\pi}e^{\frac{b^{2}}{4a}}$ \cite{Erf(ifny)}. 

\textbf{Theorem 2}

\begin{equation}
\begin{array}{ll}
V\left(0,\,2,\, a,\, b,\, x\right)=\sqrt{\frac{2}{\pi}}\int e^{-ax^{2}}b^{n+\frac{1}{2}}x^{0+\frac{1}{2}}\frac{1}{2}G_{0,2}^{2,0}\left(\frac{b^{2}x^{2}}{4}|\begin{array}{c}
\frac{5}{4},-\frac{5}{4}\end{array}\right)\, dx\\
=\sqrt{\frac{2}{\pi}}\int e^{-ax^{2}}b^{2+\frac{1}{2}}x^{0+\frac{1}{2}}K_{\frac{5}{2}}(xb)\, dx\\
=b^{2}\int\left(\frac{3}{x^{2}b^{2}}+\frac{3}{xb}+1\right)e^{-xb-ax^{2}}\, dx\\
=-\frac{3e^{-ax^{2}-bx}}{x}-\frac{\sqrt{\pi}\left(6a-b^{2}\right)e^{\frac{b^{2}}{4a}}\text{erf}\left(\frac{2ax+b}{2\sqrt{a}}\right)}{2\sqrt{a}}\end{array}\:\:.\label{eq:V(0,2,x) Theorem 2}\end{equation}

\textbf{Proof:}

In this second case one again has $\upsilon=2$, but with an additional
constant term and a coefficient of 3 shared by the first two terms.
Using \cite{PBM1} (p. 140 No. 1.3.3.17)

\begin{equation}
\int e^{-ax^{2}-bx}\, dx=\frac{\sqrt{\pi}}{2\sqrt{a}}e^{\frac{b^{2}}{4a}}\text{erf}\left(\frac{2ax+b}{2\sqrt{a}}\right)\label{eq:PBM1 p. 140 No. 1.3.3.17}\end{equation}
for the constant term, one has

\begin{equation}
\begin{array}{ll}
V\left(0,\,2,\, a,\, b,\, x\right)=\int e^{-bx-ax^{2}}\left(\frac{3}{x^{2}}+\frac{3b}{x}+b^{2}\right)\, dx=-3\frac{e^{-ax^{2}-\text{bx}}}{x}-\int\left(6a-b^{2}\right)e^{-ax^{2}-bx}\, dx-\int(\frac{3b}{x}-\frac{3b}{x})e^{-ax^{2}-bx}\, dx\\
=-\frac{3e^{-ax^{2}-bx}}{x}-\frac{\sqrt{\pi}\left(6a-b^{2}\right)e^{\frac{b^{2}}{4a}}\text{erf}\left(\frac{2ax+b}{2\sqrt{a}}\right)}{2\sqrt{a}}\end{array}\:,\label{eq:V(0,2,x) proof}\end{equation}
thus completing the proof. $\square$

\textbf{Theorem 3}

\begin{equation}
\begin{array}{l}
V\left(1,\,3,\, a,\, b,\, x\right)=\sqrt{\frac{2}{\pi}}\int e^{-ax^{2}}b^{3+\frac{1}{2}}x^{1+\frac{1}{2}}\frac{1}{2}G_{0,2}^{2,0}\left(\frac{x^{2}b^{2}}{4}|\begin{array}{c}
\frac{7}{4},-\frac{7}{4}\end{array}\right)\, dx\\
=\sqrt{\frac{2}{\pi}}\int e^{-ax^{2}}b^{3+\frac{1}{2}}x^{1+\frac{1}{2}}K_{\frac{7}{2}}(xb)\, dx\\
=b^{3}\int x\left(\frac{15}{x^{3}b^{3}}+\frac{15}{x^{2}b^{2}}+\frac{6}{xb}+1\right)e^{-xb-ax^{2}}\, dx\:\\
=-e^{-ax^{2}-bx}\left(\frac{15}{x}+\frac{b^{3}}{2a}\right)-\frac{\sqrt{\pi}\left(60a^{2}-12ab^{2}+b^{4}\right)e^{\frac{b^{2}}{4a}}\text{erf}\left(\frac{2ax+b}{2\sqrt{a}}\right)}{4a^{3/2}}\end{array}\:,\label{eq:V(1,3,x) Theorem 3}\end{equation}

\textbf{Proof:}

One again has $\upsilon=2$, and an additional constant term along
with a shared coefficient of 15 for the first two terms. In this case
there is also a term with positive power of x. The latter, and higher
powers, may be found by combining four integrals in Prudnikov, Brychkov,
and Marichev~ \cite{PBM1} (p. 139-40 No. 1.3.3.1, 5, 6, and 19),
but care must be taken since No. 19 has the opposite sign for \emph{b}
from what we have been using\emph{,} and \emph{a} is squared in three
of the four integrals. For this reason I have provided a Mathematica
program with consistent parameters in Appendix A based on Korotkov's
more straightforward general-power result \cite{Korotkov} (p. 4 No.
1.2.4 with $a\rightarrow\alpha^{2}$ and $b\rightarrow-\beta$)) .
One may also take a series of derivatives of (\ref{eq:PBM1 p. 140 No. 1.3.3.17})
with respect to \emph{b}, the first three of which are

\begin{equation}
\int xe^{-ax^{2}-bx}\, dx=-\frac{\partial}{\partial b}\int e^{-ax^{2}-bx}\, dx=-\frac{\partial}{\partial b}\frac{\sqrt{\pi}}{2\sqrt{a}}e^{\frac{b^{2}}{4a}}\text{erf}\left(\frac{2ax+b}{2\sqrt{a}}\right)=-\frac{e^{-ax^{2}-bx}}{2a}-\frac{\sqrt{\pi}b}{4a^{3/2}}e^{\frac{b^{2}}{4a}}\text{erf}\left(\frac{2ax+b}{2\sqrt{a}}\right)\:,\label{eq:int x exp}\end{equation}

\begin{equation}
\int x^{2}e^{-ax^{2}-bx}\, dx=-e^{-ax^{2}-bx}\left(\frac{x}{2a}-\frac{b}{4a^{2}}\right)+\frac{\sqrt{\pi}\left(2a+b^{2}\right)e^{\frac{b^{2}}{4a}}\text{erf}\left(\frac{2ax+b}{2\sqrt{a}}\right)}{8a^{5/2}}\label{eq:int x^2 exp}\end{equation}
(both of which are given in Korotkov \cite{Korotkov} (p. 3 No. 1.2.2,3
with $a\rightarrow\alpha^{2}$ and $b\rightarrow-\beta$)) and

\begin{equation}
\begin{array}{l}
\int x^{3}e^{-ax^{2}-bx}\, dx=-\frac{\sqrt{\pi}b\left(6a+b^{2}\right)e^{\frac{b^{2}}{4a}}\text{erf}\left(\frac{2ax+b}{2\sqrt{a}}\right)}{16a^{7/2}}-e^{-ax^{2}-bx}\left(\frac{b^{2}}{8a^{3}}-\frac{bx}{4a^{2}}+\frac{1}{2a^{2}}+\frac{x^{2}}{2a}\right)\end{array}\:.\label{eq:int x^3 exp}\end{equation}

Using (\ref{eq:int x exp}) one finds that

\begin{equation}
\begin{array}{ll}
V\left(1,\,3,\, a,\, b,\, x\right)=\int\left(\frac{15}{x^{2}}+\frac{15b}{x}+6b^{2}+b^{3}x\right)e^{-ax^{2}-bx}\, dx=-15\frac{e^{-ax^{2}-\text{bx}}}{x}-\int\left(30a-6b^{2}\right)e^{-ax^{2}-bx}\, dx-\int(\frac{15b}{x}-\frac{15b}{x})e^{-ax^{2}-bx}\, dx\\
\qquad\qquad\qquad+b^{3}\int xe^{-ax^{2}-bx}\, dx\\
=-\frac{15e^{-ax^{2}-bx}}{x}-\frac{\sqrt{\pi}\left(30a-6b^{2}\right)}{2\sqrt{a}}e^{\frac{b^{2}}{4a}}\text{erf}\left(\frac{2ax+b}{2\sqrt{a}}\right)+b^{3}\left(-\frac{\sqrt{\pi}be^{\frac{b^{2}}{4a}}\text{erf}\left(\frac{2ax+b}{2\sqrt{a}}\right)}{4a^{3/2}}-\frac{e^{-ax^{2}-bx}}{2a}\right)\end{array}\:,\label{eq:V(1,3,x) proof}\end{equation}
and rearrangement completes the proof. $\square$

\textbf{Theorem 4}

\begin{equation}
\begin{array}{ll}
V\left(2,\,4,\, a,\, b,\, x\right)=\sqrt{\frac{2}{\pi}}\int e^{-ax^{2}}b^{4+\frac{1}{2}}x^{2+\frac{1}{2}}\frac{1}{2}G_{0,2}^{2,0}\left(\frac{x^{2}b^{2}}{4}|\begin{array}{c}
\frac{9}{4},-\frac{9}{4}\end{array}\right)\, dx\\
=\sqrt{\frac{2}{\pi}}\int e^{-ax^{2}}b^{4+\frac{1}{2}}x^{2+\frac{1}{2}}K_{\frac{9}{2}}(xb)\,\, dx\:\\
=b^{4}\int x^{2}\left(\frac{105}{x^{4}b^{4}}+\frac{105}{x^{3}b^{3}}+\frac{45}{x^{2}b^{2}}+\frac{10}{xb}+1\right)e^{-ax^{2}-bx}\, dx\\
=-e^{-ax^{2}-bx}\left(\frac{20ab^{3}-b^{5}}{4a^{2}}+\frac{b^{4}x}{2a}+\frac{105}{x}\right)\\
-\frac{\sqrt{\pi}\left(840a^{3}-180a^{2}b^{2}+18ab^{4}-b^{6}\right)}{8a^{5/2}}e^{\frac{b^{2}}{4a}}\text{erf}\left(\frac{2ax+b}{2\sqrt{a}}\right)\end{array}\:,\label{eq:V(2,4,x) Theorem 4}\end{equation}

\textbf{Proof:}

Once again has $\upsilon=2$ so that one can use (\ref{eq:uv-intvdu})
for the first two terms, (\ref{eq:PBM1 p. 140 No. 1.3.3.17}) for
the constant term and (\ref{eq:int x exp}) and (\ref{eq:int x^2 exp}),
resp., for the final two terms having powers of x. One finds that

\begin{equation}
\begin{array}{ll}
V\left(2,\,4\, a,\, b,\, x\right)=\int\left(\frac{105}{x^{2}}+\frac{105b}{x}+45b^{2}+10b^{3}x+b^{4}x^{2}\right)e^{-ax^{2}-bx}\, dx\\
=-105\frac{e^{-ax^{2}-\text{bx}}}{x}-\int\left(210a-45b^{2}\right)e^{-ax^{2}-bx}\, dx-\int(\frac{105b}{x}-\frac{105b}{x})e^{-ax^{2}-bx}\, dx\\
\qquad\qquad\qquad+10b^{3}\int xe^{-ax^{2}-bx}\, dx+b^{4}\int x^{2}e^{-ax^{2}-bx}\, dx\\
=-\frac{105e^{-ax^{2}-bx}}{x}-\frac{\sqrt{\pi}\left(210a-45b^{2}\right)e^{\frac{b^{2}}{4a}}\text{erf}\left(\frac{2ax+b}{2\sqrt{a}}\right)}{2\sqrt{a}}+10b^{3}\left(-\frac{\sqrt{\pi}be^{\frac{b^{2}}{4a}}\text{erf}\left(\frac{2ax+b}{2\sqrt{a}}\right)}{4a^{3/2}}-\frac{e^{-ax^{2}-bx}}{2a}\right)\\
+b^{4}\left(\frac{\sqrt{\pi}\left(2a+b^{2}\right)e^{\frac{b^{2}}{4a}}\text{erf}\left(\frac{2ax+b}{2\sqrt{a}}\right)}{8a^{5/2}}+e^{-ax^{2}-bx}\left(\frac{b}{4a^{2}}-\frac{x}{2a}\right)\right)\end{array}\:,\label{eq:V(2,4,x) proof}\end{equation}
and rearrangement completes the proof. $\square$

The final integral we shall do in this series is

\textbf{Theorem 5}

\begin{equation}
\begin{array}{ll}
V\left(3,\,5,\, a,\, b,\, x\right)=\sqrt{\frac{2}{\pi}}\int e^{-ax^{2}}b^{5+\frac{1}{2}}x^{3+\frac{1}{2}}\frac{1}{2}G_{0,2}^{2,0}\left(\frac{x^{2}b^{2}}{4}|\begin{array}{c}
\frac{11}{4},-\frac{11}{4}\end{array}\right)\, dx\\
=\sqrt{\frac{2}{\pi}}\int e^{-ax^{2}}b^{5+\frac{1}{2}}x^{3+\frac{1}{2}}K_{\frac{11}{2}}(xb)\,\, dx\\
=b^{5}\int x^{3}\left(\frac{945}{x^{5}b^{5}}+\frac{945}{x^{4}b^{4}}+\frac{420}{x^{3}b^{3}}+\frac{105}{x^{2}b^{2}}+\frac{15}{xb}+1\right)e^{-ax^{2}-bx}\,\, dx\\
= \hspace{0.2 cm} e^{-ax^{2}-bx}\left(-\frac{b^{7}}{8a^{3}}+\frac{b^{6}x}{4a^{2}}+\frac{13b^{5}}{4a^{2}}-\frac{b^{5}x^{2}}{2a}-\frac{15b^{4}x}{2a}-\frac{105b^{3}}{2a}-\frac{945}{x}\right)\\
-\frac{\sqrt{\pi}\left(15120a^{4}-3360a^{3}b^{2}+360a^{2}b^{4}-24ab^{6}+b^{8}\right)e^{\frac{b^{2}}{4a}}\text{erf}\left(\frac{2ax+b}{2\sqrt{a}}\right)}{16a^{7/2}}\end{array}\:,\label{eq:V(3,5,x) Theorem 5}\end{equation}

\textbf{Proof:}

The only substantive change from the Theorem 4 proof is the addition
of a third power of x for which one needs (\ref{eq:int x^3 exp}).
Then

\begin{equation}
\begin{array}{ll}
V\left(3,\,5,\, a,\, b,\, x\right)=\int\left(\frac{945}{x^{2}}+\frac{945b}{x}+420b^{2}+105b^{3}x+15b^{4}x^{2}+b^{5}x^{3}\right)e^{-ax^{2}-bx}\, dx\\
=-945\frac{e^{-ax^{2}-\text{bx}}}{x}-\int\left(1890a-420b^{2}\right)e^{-ax^{2}-bx}\, dx\\
-\int(\frac{945b}{x}-\frac{945b}{x})e^{-ax^{2}-bx}\, dx+105b^{3}\int xe^{-ax^{2}-bx}\, dx+15b^{4}\int x^{2}e^{-ax^{2}-bx}\, dx+b^{5}\int x^{3}e^{-ax^{2}-bx}\, dx\\
=-\frac{945e^{-ax^{2}-bx}}{x}-\frac{\sqrt{\pi}\left(1890a-420b^{2}\right)e^{\frac{b^{2}}{4a}}\text{erf}\left(\frac{2ax+b}{2\sqrt{a}}\right)}{2\sqrt{a}}+105b^{3}\left(-\frac{\sqrt{\pi}be^{\frac{b^{2}}{4a}}\text{erf}\left(\frac{2ax+b}{2\sqrt{a}}\right)}{4a^{3/2}}-\frac{e^{-ax^{2}-bx}}{2a}\right)\\
+15b^{4}\left(\frac{\sqrt{\pi}\left(2a+b^{2}\right)e^{\frac{b^{2}}{4a}}\text{erf}\left(\frac{2ax+b}{2\sqrt{a}}\right)}{8a^{5/2}}+e^{-ax^{2}-bx}\left(\frac{b}{4a^{2}}-\frac{x}{2a}\right)\right)\\
+b^{5}\left(e^{-ax^{2}-bx}\left(-\frac{b^{2}}{8a^{3}}+\frac{bx}{4a^{2}}-\frac{1}{2a^{2}}-\frac{x^{2}}{2a}\right)-\frac{\sqrt{\pi}b\left(6a+b^{2}\right)e^{\frac{b^{2}}{4a}}\text{erf}\left(\frac{2ax+b}{2\sqrt{a}}\right)}{16a^{7/2}}\right)\end{array}\:,\label{eq:V(3,5,x) proof}\end{equation}
and rearrangement completes the proof. $\square$

\section{Integrals with $h=n-4$ }

One may also investigate indefinite
integrals of Macdonald functions multiplied by powers such that the
largest inverse power of \emph{x} is 4, such as

\textbf{Theorem 6}

\begin{equation}
\begin{array}{ll}
V\left(-2,\,2,\, a,\, b,\, x\right)=\sqrt{\frac{2}{\pi}}\int e^{-ax^{2}}b^{2+\frac{1}{2}}x^{-2+\frac{1}{2}}\frac{1}{2}G_{0,2}^{2,0}\left(\frac{b^{2}x^{2}}{4}|\begin{array}{c}
\frac{5}{4},-\frac{5}{4}\end{array}\right)\, dx\\
=\sqrt{\frac{2}{\pi}}\int e^{-ax^{2}}b^{2+\frac{1}{2}}x^{-2+\frac{1}{2}}K_{\frac{5}{2}}(xb)\,\, dx\\
=b^{2}\int\frac{1}{x^{2}}\left(\frac{3}{x^{2}b^{2}}+\frac{3}{xb}+1\right)e^{-xb-ax^{2}}\, dx\\
=e^{-ax^{2}-bx}\left(\frac{2a}{x}-\frac{b}{x^{2}}-\frac{1}{x^{3}}\right)+2\sqrt{\pi}a^{3/2}e^{\frac{b^{2}}{4a}}\text{erf}\left(\frac{2ax+b}{2\sqrt{a}}\right)\end{array}\:\label{eq:V(-2,2,x) Theorem 6}\end{equation}

\textbf{Proof:}

Unfortunately, integration by parts for the first term of \begin{equation}
\int\left(\frac{3}{x^{4}}+\frac{3b}{x^{3}}+\frac{b^{2}}{x^{2}}\right)e^{-ax^{2}-bx}\, dx\label{eq:x^-4 integral}\end{equation}
 using (\ref{eq:uv-intvdu}) with $\nu=4$ yields

\begin{equation}
-\frac{e^{-ax^{2}-\text{bx}}}{x^{3}}-\int\left(\frac{2ae^{-ax^{2}-bx}}{x^{2}}+\frac{be^{-ax^{2}-bx}}{x^{3}}\right)\, dx\:,\label{eq:uv-intvdu w nu=00003D4}\end{equation}
whose second term now has the wrong coefficient to cancel the second
term of (\ref{eq:x^-4 integral}).

However, noticing that minus the indefinite integral of the third
line of (\ref{eq:V(0,2,x) Theorem 2}) with respect to \emph{a} yields
(\ref{eq:x^-4 integral}), one might hope one could also find minus
the indefinite integral of last line with respect to \emph{a} in order
to prove Theorem 6. This is a high bar since indefinite integrals
of error functions with complicated arguments are as scarce as indefinite
integrals of Macdonald functions. There are none in section 4.1 of
one of the most conprehensive set of tables given by Ng and Geller
\cite{Ng and Geller}, nor are the correct combinations to be found
in sections 1.3-8 of Korotkov \cite{Korotkov}, the other comprehensive
table. 

One can, however, write the error function in an integral representation
\cite{PBM1} (p. 140 No. 1.3.3.17)

\begin{equation}
\text{erf}\left(\frac{2ax+b}{2\sqrt{a}}\right)=\frac{2\sqrt{a}e^{-\frac{b^{2}}{4a}}}{\sqrt{\pi}}\int e^{-ax^{2}-bx}\, dx\label{eq:erf as int rep}\end{equation}
so that the second term on the last line of (\ref{eq:V(0,2,x) Theorem 2})
can be respresented as\begin{align}
-\frac{\left(\sqrt{\pi}\left(6a-b^{2}\right)e^{\frac{b^{2}}{4a}}\right)\text{erf}\left(\frac{2ax+b}{2\sqrt{a}}\right)}{2\sqrt{a}} & =-\frac{\left(\sqrt{\pi}\left(6a-b^{2}\right)e^{\frac{b^{2}}{4a}}\right)}{2\sqrt{a}}\left(\frac{2\sqrt{a}e^{-\frac{b^{2}}{4a}}}{\sqrt{\pi}}\int e^{-ax^{2}-bx}\, dx\right)\label{eq:V02 erf term as int rep}\\
= & -\left(6a-b^{2}\right)\int e^{-ax^{2}-bx}\, dx\nonumber \end{align}

One can easily find the negative of the indefinite integral of this
integrand with respect to \emph{a} \cite{GR5}(p. 112 Nos. 2.311 and
2.322.1) \begin{equation}
\left(-1\right)^{2}\int\left(6a-b^{2}\right)e^{-ax^{2}-bx}\, da=e^{-ax^{2}-bx}\left(-\frac{6a}{x^{2}}+\frac{b^{2}}{x^{2}}-\frac{6}{x^{4}}\right)\label{eq:a-int of V02 erf term}\end{equation}
but the integral of this with respect to \emph{x,} to remove the introduced
integral representation, takes some work. A three-step integration
by parts for this integral is given in Appendix B, along with other
auxiliary indefinite integrals.  It gives

\begin{align}
\left(-1\right)^{2}\int\frac{\left(\sqrt{\pi}\left(6a-b^{2}\right)e^{\frac{b^{2}}{4a}}\right)\text{erf}\left(\frac{2ax+b}{2\sqrt{a}}\right)}{2\sqrt{a}}\, da & =\int e^{-ax^{2}-bx}\left(-\frac{6a}{x^{2}}+\frac{b^{2}}{x^{2}}-\frac{6}{x^{4}}\right)\, dx\label{eq:x-int of V02 erf term}\\
 & =2\sqrt{\pi}a^{3/2}e^{\frac{b^{2}}{4a}}\text{erf}\left(\frac{2ax+b}{2\sqrt{a}}\right)+e^{-ax^{2}-bx}\left(\frac{2a}{x}-\frac{b}{x^{2}}+\frac{2}{x^{3}}\right)\nonumber \end{align}

The the negative of the indefinite integral with respect to \emph{a}
of the first term on the last line of (\ref{eq:V(0,2,x) Theorem 2})
is more easily done:\begin{equation}
\left(-1\right)^{2}\int\frac{3e^{-ax^{2}-bx}}{x}\, da=-\frac{3e^{-ax^{2}-bx}}{x^{3}}\:.\label{eq:int wrt a of 1st term of V02}\end{equation}
When this is added to the final term of (\ref{eq:x-int of V02 erf term}),
it gives (\ref{eq:V(-2,2,x) Theorem 6}), thus completing the proof.
$\square$

\textbf{Theorem 7}

\begin{equation}
\begin{array}{l}
V\left(-1,\,3,\, a,\, b,\, x\right)=\sqrt{\frac{2}{\pi}}\int e^{-ax^{2}}b^{3+\frac{1}{2}}x^{-1+\frac{1}{2}}\frac{1}{2}G_{0,2}^{2,0}\left(\frac{x^{2}b^{2}}{4}|\begin{array}{c}
\frac{7}{4},-\frac{7}{4}\end{array}\right)\, dx\\
=\sqrt{\frac{2}{\pi}}\int e^{-ax^{2}}b^{-1+\frac{1}{2}}x^{-1+\frac{1}{2}}K_{\frac{7}{2}}(xb)\, dx\\
=b^{3}\int\frac{1}{x}\left(\frac{15}{x^{3}b^{3}}+\frac{15}{x^{2}b^{2}}+\frac{6}{xb}+1\right)e^{-ax^{2}-bx}\, dx\\
=e^{-ax^{2}-bx}\left(\frac{10a}{x}-\frac{b^{2}}{x}-\frac{5b}{x^{2}}-\frac{5}{x^{3}}\right)+\sqrt{\pi}\sqrt{a}\left(10a-b^{2}\right)e^{\frac{b^{2}}{4a}}\text{erf}\left(\frac{2ax+b}{2\sqrt{a}}\right)\end{array}\:,\label{eq:V(-1,3,x) Theorem 7}\end{equation}

\textbf{Proof:}

Since we know that the \emph{x}-dependence of the final form will
be similar to the final line of Theorem 6, the easiest path to a proof
is to take the derivative of some posited form of the result with
respect to \emph{x,}

\begin{align}
 & \frac{\partial}{\partial x}\left(e^{-ax^{2}-bx}\left(\frac{c}{x^{3}}+\frac{d}{x^{2}}+\frac{f}{x}\right)+\sqrt{\pi}ge^{\frac{b^{2}}{4a}}\text{erf}\left(\frac{2ax+b}{2\sqrt{a}}\right)\right)\:,\label{eq:test D or x^-3 through x^-1 +erf}\\
= & e^{-ax^{2}-bx}\left(-\frac{3c}{x^{4}}+\frac{-bc-2d}{x^{3}}+\frac{-2ac-bd-f}{x^{2}}+\frac{-2ad-bf}{x}-2af+2\sqrt{a}g\right)\nonumber \end{align}
and match coefficients with the integrand in the third line of Theorem
7. Doing the identification from the left-most term to right-most
in the third line of Theorem 7, one sees that the set $\left\{ c\to-5,\, d\to-5b, \right.$ $\left.\, f\to10a-b^{2}, \right.$ $\left.\, g\to\sqrt{a}\left(10a-b^{2}\right)\right\} $
is the correct solution to complete the proof. $\square$

\textbf{Theorem 8}

\begin{equation}
\begin{array}{ll}
V\left(0,\,4,\, a,\, b,\, x\right)=\sqrt{\frac{2}{\pi}}\int e^{-ax^{2}}b^{4+\frac{1}{2}}x^{0+\frac{1}{2}}\frac{1}{2}G_{0,2}^{2,0}\left(\frac{x^{2}b^{2}}{4}|\begin{array}{c}
\frac{9}{4},-\frac{9}{4}\end{array}\right)\, dx\\
=\sqrt{\frac{2}{\pi}}\int e^{-ax^{2}}b^{4+\frac{1}{2}}x^{0+\frac{1}{2}}K_{\frac{9}{2}}(xb)\,\, dx\:\\
=b^{4}\int\left(\frac{105}{x^{4}b^{4}}+\frac{105}{x^{3}b^{3}}+\frac{45}{x^{2}b^{2}}+\frac{10}{xb}+1\right)e^{-ax^{2}-bx}\,\, dx\\
=e^{-ax^{2}-bx}\left(-\frac{35}{x^{3}}-\frac{35b}{x^{2}}+\frac{10\left(7a-b^{2}\right)}{x}\right)\\
+\frac{\sqrt{\pi}\left(140a^{2}-20ab^{2}+b^{4}\right)}{2\sqrt{a}}e^{\frac{b^{2}}{4a}}\text{erf}\left(\frac{2ax+b}{2\sqrt{a}}\right)\end{array}\:,\label{eq:V(0,4,x) Theorem 8}\end{equation}

\textbf{Proof:}

One again can match the coeffieicnts in (\ref{eq:test D or x^-3 through x^-1 +erf})
to the third line of Theorem 8, taking care to account for its new
constant term. Doing the identification from the left-most term to
right-most, one sees that the set $\left\{ c\to-35,\, d\to-35b, \, f\to10\left(7a-b^{2}\right),\, g\to\frac{140a^{2}-20ab^{2}+b^{4}}{2\sqrt{a}}\right\} $
is the correct solution to complete the proof, with the final term
within \emph{g} providing the new constant term in the third line
of Theorem 8. $\square$

In the prior work \cite{stra26a}, the Macdonald functions for both
$n=4$ and $5$ were mistakenly labeled with index $K_{\frac{7}{2}}(xb)$.

\textbf{Theorem 9}

\begin{equation}
\begin{array}{ll}
V\left(1,\,5,\, a,\, b,\, x\right)=\sqrt{\frac{2}{\pi}}\int e^{-ax^{2}}b^{5+\frac{1}{2}}x^{1+\frac{1}{2}}\frac{1}{2}G_{0,2}^{2,0}\left(\frac{x^{2}b^{2}}{4}|\begin{array}{c}
\frac{11}{4},-\frac{11}{4}\end{array}\right)\, dx\\
=\sqrt{\frac{2}{\pi}}\int e^{-ax^{2}}b^{5+\frac{1}{2}}x^{1+\frac{1}{2}}K_{\frac{11}{2}}(xb)\,\, dx\\
=b^{5}\int x\left(\frac{945}{x^{5}b^{5}}+\frac{945}{x^{4}b^{4}}+\frac{420}{x^{3}b^{3}}+\frac{105}{x^{2}b^{2}}+\frac{15}{xb}+1\right)e^{-ax^{2}-bx}\,\, dx\\
=e^{-ax^{2}-bx}\left(-\frac{b^{5}}{2a}+\frac{105\left(6a-b^{2}\right)}{x}-\frac{315b}{x^{2}}-\frac{315}{x^{3}}\right)\\
+\frac{\sqrt{\pi}\left(2520a^{3}-420a^{2}b^{2}+30ab^{4}-b^{6}\right)}{4a^{3/2}}e^{\frac{b^{2}}{4a}}\text{erf}\left(\frac{2ax+b}{2\sqrt{a}}\right)\end{array}\:,\label{eq:V(1,5,x) Theorem 9}\end{equation}

\textbf{Proof:}

In the present case, the third line of Theorem 9 has an $x^{1}$ term
that requires a modificantion of (\ref{eq:test D or x^-3 through x^-1 +erf})
to include a constant term \emph{h}:

\begin{align}
 & \frac{\partial}{\partial x}\left(e^{-ax^{2}-bx}\left(\frac{c}{x^{3}}+\frac{d}{x^{2}}+\frac{f}{x}+h\right)+\sqrt{\pi}ge^{\frac{b^{2}}{4a}}\text{erf}\left(\frac{2ax+b}{2\sqrt{a}}\right)\right)\label{eq:test D of x^-3 through x^0 +erf}\\
= & e^{-ax^{2}-bx}\left(\frac{945}{x^{4}}+\frac{945b}{x^{3}}+\frac{630a+315b^{2}-f}{x^{2}}+\frac{630ab-bf}{x}-2af+2\sqrt{a}g-2ahx-bh\right)\nonumber \end{align}
but the identification procedure remains the same. One sees that the
set $\left\{ c\to-315,\, d\to-315b, \right.$ 

\noindent
$\left. \, f\to105\left(6a-b^{2}\right),\, g\to\frac{1260a^{2}-210ab^{2}+15b^{4}+bh}{2\sqrt{a}}\right\} $
followed by $h\to-\frac{b^{5}}{2a}$ , completes the proof. $\square$

\section{Integrals with $h=n-6$ }

One may also investigate indefinite
integrals of Macdonald functions multiplied by powers such that the
largest inverse power of \emph{x} is 6, such as

\textbf{Theorem 10}

\begin{equation}
\begin{array}{l}
V\left(-3,\,3,\, a,\, b,\, x\right)=\sqrt{\frac{2}{\pi}}\int e^{-ax^{2}}b^{3+\frac{1}{2}}x^{-3+\frac{1}{2}}\frac{1}{2}G_{0,2}^{2,0}\left(\frac{x^{2}b^{2}}{4}|\begin{array}{c}
\frac{7}{4},-\frac{7}{4}\end{array}\right)\, dx\\
=\sqrt{\frac{2}{\pi}}\int e^{-ax^{2}}b^{3+\frac{1}{2}}x^{-3+\frac{1}{2}}K_{\frac{7}{2}}(xb)\,\, dx\\
=b^{3}\int\frac{1}{x^{3}}\left(\frac{15}{x^{3}b^{3}}+\frac{15}{x^{2}b^{2}}+\frac{6}{xb}+1\right)e^{-ax^{2}-bx}\, dx\\
=e^{-ax^{2}-bx}\left(-\frac{3}{x^{5}}-\frac{3b}{x^{4}}+\frac{2a-b^{2}}{x^{3}}+\frac{2ab}{x^{2}}-\frac{4a^{2}}{x}\right)-4\sqrt{\pi}a^{5/2}e^{\frac{b^{2}}{4a}}\text{erf}\left(\frac{2ax+b}{2\sqrt{a}}\right)\end{array}\:.\label{eq:V(-3,3,x) Theorem 10}\end{equation}

\textbf{Proof:}

Noticing that minus the indefinite integral of the third line of (\ref{eq:V(-1,3,x) Theorem 7})
with respect to \emph{a} yields the third line of (\ref{eq:V(-3,3,x) Theorem 10}),
one hopes to find the indefinite integral of last line of (\ref{eq:V(-1,3,x) Theorem 7})
with respect to \emph{a} in order to prove Theorem 10. Using the integral
representation (\ref{eq:erf as int rep}), the second term on the
last line of (\ref{eq:V(-1,3,x) Theorem 7}) can be respresented as\begin{align}
\sqrt{\pi}\sqrt{a}\left(10a-b^{2}\right)e^{\frac{b^{2}}{4a}}\text{erf}\left(\frac{2ax+b}{2\sqrt{a}}\right)=2a\left(10a-b^{2}\right)\int e^{-ax^{2}-bx}\, dx\label{eq:V(-1,3) erf term as int rep}\end{align}

One can easily find the negative of indefinite integral of this integrand
with respect to \emph{a} \cite{GR5}(p. 112 Nos. 2.322.1,2) \begin{equation}
-\int2a\left(10a-b^{2}\right)e^{-ax^{2}-bx}\, da=e^{-ax^{2}-bx}\left(\frac{40}{x^{6}}+\frac{40a-2b^{2}}{x^{4}}+\frac{20a^{2}-2ab^{2}}{x^{2}}\right)\label{eq:a-int of V(-1,3) erf term}\end{equation}
but the integral of this with respect to \emph{x,} to remove the introduced
integral representation, takes extensive work. A five-step integration
by parts for this integral is given in Appendix B, such that

\begin{align}
& \left(-1\right)\int\sqrt{\pi}\sqrt{a}\left(10a-b^{2}\right)e^{\frac{b^{2}}{4a}}\text{erf}\left(\frac{2ax+b}{2\sqrt{a}}\right)\, da \nonumber \\
& =\int e^{-ax^{2}-bx}\left(\frac{40}{x^{6}}+\frac{40a-2b^{2}}{x^{4}}+\frac{20a^{2}-2ab^{2}}{x^{2}}\right)\, dx\:.\label{eq:x-int of V(-1,3) erf term}\\
 & =e^{-ax^{2}-bx}\left(-\frac{4a^{2}}{x}+\frac{2ab}{x^{2}}-\frac{8a}{x^{3}}+\frac{2b}{x^{4}}-\frac{8}{x^{5}}\right)-4\sqrt{\pi}a^{5/2}e^{\frac{b^{2}}{4a}}\text{erf}\left(\frac{2ax+b}{2\sqrt{a}}\right)\nonumber \end{align}

The negative of indefinite integral of of the first term on the last
line of (\ref{eq:V(-1,3,x) Theorem 7}) with respect to \emph{a} is
more easily done \cite{GR5} (p. 112 Nos. 2.311 and 2.322.1):\begin{equation}
-\int e^{-ax^{2}-bx}\left(\frac{10a}{x}-\frac{b^{2}}{x}-\frac{5b}{x^{2}}-\frac{5}{x^{3}}\right)\, da=e^{-ax^{2}-bx}\left(\frac{10a}{x^{3}}-\frac{b^{2}}{x^{3}}-\frac{5b}{x^{4}}+\frac{5}{x^{5}}\right)\:.\label{eq:int wrt a of 1st term of V(-1,3)}\end{equation}
When this is added to the first term of the last line of (\ref{eq:x-int of V(-1,3) erf term}),
it gives the last line of (\ref{eq:V(-3,3,x) Theorem 10}), thus completing
the proof. $\square$

\textbf{Theorem 11}

\begin{equation}
\begin{array}{ll}
V\left(-2,\,4,\, a,\, b,\, x\right)=\sqrt{\frac{2}{\pi}}\int e^{-ax^{2}}b^{4+\frac{1}{2}}x^{-2+\frac{1}{2}}\frac{1}{2}G_{0,2}^{2,0}\left(\frac{x^{2}b^{2}}{4}|\begin{array}{c}
\frac{9}{4},-\frac{9}{4}\end{array}\right)\, dx\\
=\sqrt{\frac{2}{\pi}}\int e^{-ax^{2}}b^{4+\frac{1}{2}}x^{-2+\frac{1}{2}}K_{\frac{9}{2}}(xb)\,\, dx\:\\
=b^{4}\int\frac{1}{x^{2}}\left(\frac{105}{x^{4}b^{4}}+\frac{105}{x^{3}b^{3}}+\frac{45}{x^{2}b^{2}}+\frac{10}{xb}+1\right)e^{-ax^{2}-bx}\,\, dx\\
=\left(e^{-ax^{2}-bx}\left(\frac{2ab^{2}-28a^{2}}{x}+\frac{14ab-b^{3}}{x^{2}}+\frac{14a-8b^{2}}{x^{3}}-\frac{21b}{x^{4}}-\frac{21}{x^{5}}\right)\right.\\
\left.-2\sqrt{\pi}a^{3/2}\left(14a-b^{2}\right)e^{\frac{b^{2}}{4a}}\text{erf}\left(\frac{2ax+b}{2\sqrt{a}}\right)\right)\end{array}\:,\label{eq:V(-2,4,x) Theorem 11}\end{equation}

\textbf{Proof:}

Since we know that the \emph{x}-dependence of the final form will
be similar to the final line of Theorem 10, let us take the derivative
with respect to \emph{x} of the posited form:

\begin{align}
 & \frac{\partial}{\partial x}\left(\sqrt{\pi}ge^{\frac{b^{2}}{4a}}\text{erf}\left(\frac{2ax+b}{2\sqrt{a}}\right)+e^{-ax^{2}-bx}\left(\frac{c}{x^{5}}+\frac{d}{x^{4}}+\frac{f}{x^{3}}+\frac{h}{x^{2}}+\frac{j}{x}\right)\right)\label{eq:test D or x^-6 through x^-1 +erf}\\
= & e^{-ax^{2}-bx}\left(-\frac{5c}{x^{6}}+\frac{-bc-4d}{x^{5}}+\frac{-2ac-bd-3f}{x^{4}}+\frac{-2ad-bf-2h}{x^{3}}\right.\nonumber \\
+ & \left.\frac{-2af-bh-j}{x^{2}}+\frac{-2ah-bj}{x}+2\sqrt{a}g-2aj\right)\nonumber \end{align}
and match coefficients with the third line of Theorem 11. Doing the
identification from the left-most term to right-most, one sees that
the set $\left\{ c\to-21,\, d\to-21b,\, f\to2\left(7a-4b^{2}\right),\, h\to14ab-b^{3}, \right.$ 

\noindent
$\left. \, j\to-2\left(14a^{2}-ab^{2}\right),\, g\to-2\sqrt{a}\left(14a^{2}-ab^{2}\right)\right\} $
is the correct solution to complete the proof. $\square$

\textbf{Theorem 12}

\begin{equation}
\begin{array}{ll}
V\left(-1,\,5,\, a,\, b,\, x\right)=\sqrt{\frac{2}{\pi}}\int e^{-ax^{2}}b^{5+\frac{1}{2}}x^{-1+\frac{1}{2}}\frac{1}{2}G_{0,2}^{2,0}\left(\frac{x^{2}b^{2}}{4}|\begin{array}{c}
\frac{11}{4},-\frac{11}{4}\end{array}\right)\, dx\\
=\sqrt{\frac{2}{\pi}}\int e^{-ax^{2}}b^{5+\frac{1}{2}}x^{-1+\frac{1}{2}}K_{\frac{11}{2}}(xb)\,\, dx\\
=b^{5}\int\frac{1}{x}\left(\frac{945}{x^{5}b^{5}}+\frac{945}{x^{4}b^{4}}+\frac{420}{x^{3}b^{3}}+\frac{105}{x^{2}b^{2}}+\frac{15}{xb}+1\right)e^{-ax^{2}-bx}\,\, dx\\
=e^{-ax^{2}-bx}\left(-\frac{189}{x^{5}}-\frac{189b}{x^{4}}+\frac{7\left(18a-11b^{2}\right)}{x^{3}}+\frac{14\left(9ab-b^{3}\right)}{x^{2}}+\frac{-252a^{2}+28ab^{2}-b^{4}}{x}\right)\\
-\sqrt{\pi}\sqrt{a}\left(252a^{2}-28ab^{2}+b^{4}\right)e^{\frac{b^{2}}{4a}}\text{erf}\left(\frac{2ax+b}{2\sqrt{a}}\right)\end{array}\:,\label{eq:V(-1,5,x) Theorem 12}\end{equation}

\textbf{Proof:}

Again one may take the derivative with respect to \emph{x} of the
posited form (\ref{eq:test D or x^-6 through x^-1 +erf}) and match
coefficients with the third line of Theorem 12. Doing the identification
from the left-most term to right-most, one sees that the set $\left\{ c\to-21,\, d\to-21b,\, f\to2\left(7a-4b^{2}\right),\, h\to14ab-b^{3},\, j\to-2\left(14a^{2}-ab^{2}\right),\, g\to-2\sqrt{a}\left(14a^{2}-ab^{2}\right)\right\} $
is the correct solution to complete the proof. $\square$

$ $

\section{Integrals with $h=n-8$ and \textmd{\normalsize n-10}}

One may also investigate indefinite
integrals of Macdonald functions multiplied by powers such that the
largest inverse power is 8 or 10, such as

\textbf{Theorem 13}

\begin{equation}
\begin{array}{ll}
V\left(-4,\,4,\, a,\, b,\, x\right)=\sqrt{\frac{2}{\pi}}\int e^{-ax^{2}}b^{4+\frac{1}{2}}x^{-4+\frac{1}{2}}\frac{1}{2}G_{0,2}^{2,0}\left(\frac{x^{2}b^{2}}{4}|\begin{array}{c}
\frac{9}{4},-\frac{9}{4}\end{array}\right)\, dx\\
=\sqrt{\frac{2}{\pi}}\int e^{-ax^{2}}b^{4+\frac{1}{2}}x^{-4+\frac{1}{2}}K_{\frac{9}{2}}(xb)\,\, dx\:\\
=b^{4}\int\frac{1}{x^{4}}\left(\frac{105}{x^{4}b^{4}}+\frac{105}{x^{3}b^{3}}+\frac{45}{x^{2}b^{2}}+\frac{10}{xb}+1\right)e^{-ax^{2}-bx}\,\, dx\\
=e^{-bx-ax^{2}}\left(-\frac{15}{x^{7}}-\frac{15b}{x^{6}}-\frac{6\left(-a+b^{2}\right)}{x^{5}}+\frac{6ab-b^{3}}{x^{4}}+\frac{2\left(-2a^{2}+ab^{2}\right)}{x^{3}}-\frac{4a^{2}b}{x^{2}}+\frac{8a^{3}}{x}\right)\\
+8a^{7/2}e^{\frac{b^{2}}{4a}}\sqrt{\pi}\text{erf}\left(\frac{b+2ax}{2\sqrt{a}}\right)\end{array}\:.\label{eq:V(-4,4,x) Theorem 13}\end{equation}

\textbf{Proof:}

Having now enough experience of the \emph{x}-dependence of the final
form to skip the longer proof via integration by parts, let us take
the derivative with respect to \emph{x} of an expanded posited form:

\begin{align}
 & \frac{\partial}{\partial x}\left(\sqrt{\pi}ge^{\frac{b^{2}}{4a}}\text{erf}\left(\frac{2ax+b}{2\sqrt{a}}\right)+e^{-ax^{2}-bx}\left(\frac{k}{x^{7}}+\frac{L}{x^{6}}+\frac{c}{x^{5}}+\frac{d}{x^{4}}+\frac{f}{x^{3}}+\frac{h}{x^{2}}+\frac{j}{x}\right)\right)\label{eq:test D or x^-7 through x^-1 +erf}\\
= & e^{-ax^{2}-bx}\left(-\frac{7k}{x^{8}}+\frac{-bk-6L}{x^{7}}+\frac{-2ak-bL-5c}{x^{6}}+\frac{-2aL-bc-4d}{x^{5}}\right.\nonumber \\
+ & \frac{-2ac-bd-3f}{x^{4}}+\frac{-2ad-bf-2h}{x^{3}}+\frac{-2af-bh-j}{x^{2}}+\frac{-2ah-bj}{x}\nonumber \\
+ & \left.2\sqrt{a}g-2aj\right)\nonumber \end{align}
and match coefficients with the third line of Theorem 13. Doing the
identification from the left-most term to right-most, one sees that
the set $\left\{ k\to-15,\, L\to-15b,\, c\to6\left(a-b^{2}\right),\, d\to6ab-b^{3},\, f\to-2\left(2a^{2}-ab^{2}\right), \right.$ 
$\left. \, h\to-4a^{2}b,\, j\to8a^{3},\, g\to8a^{7/2}\right\} $
is the correct solution to complete the proof. $\square$

\textbf{Theorem 14}

\begin{equation}
\begin{array}{ll}
V\left(-3,\,5,\, a,\, b,\, x\right)=\sqrt{\frac{2}{\pi}}\int e^{-ax^{2}}b^{5+\frac{1}{2}}x^{-3+\frac{1}{2}}\frac{1}{2}G_{0,2}^{2,0}\left(\frac{x^{2}b^{2}}{4}|\begin{array}{c}
\frac{11}{4},-\frac{11}{4}\end{array}\right)\, dx\\
=\sqrt{\frac{2}{\pi}}\int e^{-ax^{2}}b^{5+\frac{1}{2}}x^{-3+\frac{1}{2}}K_{\frac{11}{2}}(xb)\,\, dx\\
=b^{5}\int\frac{1}{x^{3}}\left(\frac{945}{x^{5}b^{5}}+\frac{945}{x^{4}b^{4}}+\frac{420}{x^{3}b^{3}}+\frac{105}{x^{2}b^{2}}+\frac{15}{xb}+1\right)e^{-ax^{2}-bx}\,\, dx\\
=\left(e^{-ax^{2}-bx}\left(-\frac{135}{x^{7}}-\frac{135b}{x^{6}}+\frac{54a-57b^{2}}{x^{5}}+\frac{54ab-12b^{3}}{x^{4}}\right.\right.\\
\left.+\frac{-36a^{2}+20ab^{2}-b^{4}}{x^{3}}+\frac{2ab^{3}-36a^{2}b}{x^{2}}+\frac{72a^{3}-4a^{2}b^{2}}{x}\right)\\
\left.+4\sqrt{\pi}a^{5/2}\left(18a-b^{2}\right)e^{\frac{b^{2}}{4a}}\text{erf}\left(\frac{2ax+b}{2\sqrt{a}}\right)\right)\end{array}\:,\label{eq:V(-3,5,x) Theorem 14}\end{equation}

\textbf{Proof:}

Again one may take the derivative with respect to \emph{x} of the
posited form (\ref{eq:test D or x^-7 through x^-1 +erf}). One sees
that the set $\left\{ k\to-135,\, L\to-135b,\, c\to3\left(18a-19b^{2}\right),\, d\to6\left(9ab-2b^{3}\right),\, f\to-36a^{2}+20ab^{2}-b^{4},\right.$ 

\noindent
$\left. \, h\to-2\left(18a^{2}b-ab^{3}\right),\, j\to72a^{3}-4a^{2}b^{2},\, g\to4\sqrt{a}\left(18a^{3}-a^{2}b^{2}\right)\right\} $
is the correct solution to complete the proof. $\square$

The final integral of this set has a largest inverse power of 10:

\textbf{Theorem 15}

\begin{equation}
\begin{array}{ll}
V\left(-5,\,5,\, a,\, b,\, x\right)=\sqrt{\frac{2}{\pi}}\int e^{-ax^{2}}b^{5+\frac{1}{2}}x^{-5+\frac{1}{2}}\frac{1}{2}G_{0,2}^{2,0}\left(\frac{x^{2}b^{2}}{4}|\begin{array}{c}
\frac{11}{4},-\frac{11}{4}\end{array}\right)\, dx\\
=\sqrt{\frac{2}{\pi}}\int e^{-ax^{2}}b^{5+\frac{1}{2}}x^{-5+\frac{1}{2}}K_{\frac{11}{2}}(xb)\,\, dx\\
=b^{5}\int\frac{1}{x^{5}}\left(\frac{945}{x^{5}b^{5}}+\frac{945}{x^{4}b^{4}}+\frac{420}{x^{3}b^{3}}+\frac{105}{x^{2}b^{2}}+\frac{15}{xb}+1\right)e^{-ax^{2}-bx}\,\, dx\\
=\left(\left(e^{-ax^{2}-bx}\left(-\frac{16a^{4}}{x}+\frac{8a^{3}b}{x^{2}}+\frac{2\left(ab^{3}-6a^{2}b\right)}{x^{4}}+\frac{-12a^{2}+12ab^{2}-b^{4}}{x^{5}}-\frac{4\left(a^{2}b^{2}-2a^{3}\right)}{x^{3}}\right.\right)\right.\\
\left.-\frac{10\left(b^{3}-3ab\right)}{x^{6}}-\frac{15\left(3b^{2}-2a\right)}{x^{7}}-\frac{105b}{x^{8}}-\frac{105}{x^{9}}\right)\\
\left.-16\sqrt{\pi}a^{9/2}e^{\frac{b^{2}}{4a}}\text{erf}\left(\frac{2ax+b}{2\sqrt{a}}\right)\right)\end{array}\:.\label{eq:V(-5,5,x) Theorem 15}\end{equation}

\textbf{Proof:}

One may posit an expanded form:

\begin{align}
 & \frac{\partial}{\partial x}\left(\sqrt{\pi}ge^{\frac{b^{2}}{4a}}\text{erf}\left(\frac{2ax+b}{2\sqrt{a}}\right)+e^{-ax^{2}-bx}\left(\frac{m}{x^{9}}+\frac{n}{x^{8}}+\frac{k}{x^{7}}+\frac{L}{x^{6}}+\frac{c}{x^{5}}+\frac{d}{x^{4}}+\frac{f}{x^{3}}+\frac{h}{x^{2}}+\frac{j}{x}\right)\right)\label{eq:test D or x^-9 through x^-1 +erf}\\
= & e^{-ax^{2}-bx}\left(-\frac{9m}{x^{10}}+\frac{-bm-8n}{x^{9}}+\frac{-2am-bn-7k}{x^{8}}+\frac{-2an-bk-6L}{x^{7}}\right.\nonumber \\
+ & \frac{-2ak-bL-5c}{x^{6}}+\frac{-2aL-bc-4d}{x^{5}}+\frac{-2ac-bd-3f}{x^{4}}+\frac{-2ad-bf-2h}{x^{3}}+\frac{-2af-bh-j}{x^{2}}\nonumber \\
+ & \left.\frac{-2ah-bj}{x}+2\sqrt{a}g-2aj\right)\nonumber \end{align}
and match coefficients with the third line of Theorem 15. Doing the
identification from the left-most term to right-most, one sees that
the set $\left\{ m\to-105,n\to-105b,k\to15\left(2a-3b^{2}\right),L\to10\left(3ab-b^{3}\right),\right.$ 

\noindent
$\left. c\to-12a^{2}+12ab^{2}-b^{4},d\to-2\left(6a^{2}b-ab^{3}\right),f\to4\left(2a^{3}-a^{2}b^{2}\right),h\to8a^{3}b,j\to-16a^{4},g\to-16a^{9/2}\right\} $
is the correct solution to complete the proof. $\square$

One could continue on with $n>5$, but the present set is sufficiently
comprehensive for the moment.

\section{Confirmation from computer calculus}

One always wants a check on one's results for sign errors and such,
and one might hope that a computer calculus program like \emph{Mathematica}
might provide such. Under the assumption that Mathematica would probably
convert the Macdonald function to its Meijer G-function equivalent
\cite{K_as_meijer_G} before evaluating, and in the appropriate notation
for execution, to check (\ref{eq:V(-1,1,x) Theorem 1}) I tried
\begin{quote}
Integrate{[}E\textasciicircum{}(- a (x\textasciicircum{} 2) ) x\textasciicircum{}(1/2-2
j-2 m+3 n) 1/2 MeijerG{[}\{\{\},\{\}\},\{\{1/2 (1/2+n),1/2 (-(1/2)-n)\},\{\}\},1/4
x\textasciicircum{}2 b\textasciicircum{}2{]} b\textasciicircum{}(1/2+n),x{]}/.{j->n,m->n}/.n->1 
\end{quote}
(where the notation {}``/. n->1'' directs Mathematica to make the
substitution $n=1$). The response

\begin{equation}
\sqrt{\frac{\pi}{2}}b^{3/2}\int\frac{\left(\frac{1}{\sqrt{x^{2}b^{2}}}+1\right)e^{-\sqrt{x^{2}b^{2}}-ax^{2}}}{\sqrt{x}\sqrt[4]{x^{2}b^{2}}}\, dx\:\label{eq:failed}\end{equation}
indicates the integral cannot be done (and hints that the Meijer G-function
form may not be what Mathematica uses for this integral after all).
Thinking that perhaps the factor $\sqrt{x^{2}b^{2}}$ in the exponential
and the fourth root in the denominator might be getting in the way,
I prepended two factors to remove them that each reduce to unity for
real \emph{x} and \emph{b}, 
\begin{quote}
(x\textasciicircum{}2 b\textasciicircum{}2)\textasciicircum{} (1/4)/(sqrt{[}x{]}
sqrt{[}b{]}) E\textasciicircum{}(-x b)/E\textasciicircum{}(-sqrt{[}x\textasciicircum{}2
b\textasciicircum{}2{]}) , 
\end{quote}
and this yielded the response

\begin{equation}
\sqrt{\frac{\pi}{2}}b\int\frac{\left(\frac{1}{\sqrt{x^{2}b^{2}}}+1\right)e^{-xb-ax^{2}}}{x}\, dx\:,\label{eq:failed w prepended factors}\end{equation}
 which included a term whose denominator is $\sqrt{x^{2}b^{2}}$.
So I appended the substitution command 
\begin{quote}
/. 1/sqrt{[}x\textasciicircum{}2 b\textasciicircum{}2{]} -> 1/(x b) 
\end{quote}
and both versions \emph{7} and \emph{13} of Mathematica were able
to perform the indefinite integration to give (\ref{eq:V(-1,1,x) Theorem 1}).
Neither version was able to integrate either term on the left-hand
side of the last line of (\ref{eq:V(-1,1,x) Theorem 1}) separately. 

In checking results in the prior sections for which $n>1$, one needs
to specify the substitutions 
\begin{quote}
/. {1/sqrt{[}x\textasciicircum{}2 b\textasciicircum{}2{]} -> 1/(x
b), 1/(x\textasciicircum{}2 b\textasciicircum{}2)\textasciicircum{}(3/2)
-> 1/(x\textasciicircum{}3 b\textasciicircum{}3), 1/(x\textasciicircum{}2
b\textasciicircum{}2)\textasciicircum{}(5/2) -> 1/(x\textasciicircum{}5
b\textasciicircum{}5), 1/(x\textasciicircum{}2 b\textasciicircum{}2)\textasciicircum{}(7/2)
-> 1/(x\textasciicircum{}7 b\textasciicircum{}7), 1/(x\textasciicircum{}2
b\textasciicircum{}2)\textasciicircum{}(9/2) -> 1/(x\textasciicircum{}9
b\textasciicircum{}9)} 
\end{quote}
and so on, which may be why Mathematica is unable to perform the integration
for a generic $n$. Mathematica 7 could not perform the integrals
in Section 5, but version 13 could.

I could not get \emph{Maple 2026} or \emph{Maxima} to do any such
integrals using similar notation.

A second check on displayed results for a manuscript written in \LaTeX{},
such as this one, is to copy those results from the manuscript and paste them into
Mathematica using the command 

\begin{quote}
ToExpression{[}\textquotedbl{}[paste equation
between these quotes (without the square brackets)]\textquotedbl{}, TeXForm],
\end{quote}

\noindent
As one pastes, Mathematica will ask if one wants it to insert "escapes" along with the equation material, and one responds "Yes." It will then double all of the "$\textbackslash$" characters.  E.g.,  the right-hand side of the last line of (\ref{eq:V(-1,1,x) Theorem 1}) will give

\begin{quote}
ToExpression{[}\textquotedbl{}-\textbackslash{}\textbackslash{}frac\{e\textasciicircum{}\{-ax\textasciicircum{}\{2\}-bx\}\}\{x\}-\textbackslash{}\textbackslash{}sqrt\{a\}\textbackslash{}\textbackslash{}sqrt\{\textbackslash{}\textbackslash{}pi\}   \\  e\textasciicircum{}\{\textbackslash{}\textbackslash
frac\{b\textasciicircum{}\{2\}\}\{4a\}\}\textbackslash{}\textbackslash{}text\{erf\}\textbackslash{}\textbackslash{}left(\textbackslash{}\textbackslash{}frac\{2ax+b\}\{2\textbackslash{}\textbackslash{}sqrt\{a\}\}\textbackslash{}\textbackslash{}right)\textquotedbl{}, TeXForm]
\end{quote}

\noindent
after pasting.  After executing this command, one must then convert the two lower-case \emph{e}'s to upper case \emph{E}'s (Mathematica's
short-hand for the exponential function), and in the present case
make sure to add a space between $x$ and $b$ in $E^{-xb-ax^{2}}$
so that these two are multiplied rather than treated as a single variable
named \emph{xb}. One may then take the derivative with respect to
\emph{x} to confirm a match with the integrand on the left-hand side
of the last line of (\ref{eq:V(-1,1,x) Theorem 1}). I have done so
for each of the integrals in this paper that result in analytical expressions.
Integrals involving series, below, were checked numerically.

\section{Series solution for missing integrals}

The sequence of integrals for a given value of \emph{n} in Sections
2-5 is to step \emph{h} in \emph{$x^{h}$} more negative by 2. That
is, one may find the indefinite integral $V\left(0,\,2,\, a,\, b,\, x\right)$
of Theorem 2 (\ref{eq:V(0,2,x) Theorem 2}), whose third line is $b^{2}\int\left(\frac{3}{x^{2}b^{2}}+\frac{3}{xb}+1\right)e^{-xb-ax^{2}}\, dx$,
and the next $n=2$ integral found is $V\left(-2,\,2,\, a,\, b,\, x\right)$
of Theorem 6 (\ref{eq:V(-2,2,x) Theorem 6}) whose third line is $b^{2}\int\frac{1}{x^{2}}\left(\frac{3}{x^{2}b^{2}}+\frac{3}{xb}+1\right)e^{-xb-ax^{2}}\, dx$.
Could one not then find the indefinite integral $b^{2}\int\frac{1}{x^{4}}\left(\frac{3}{x^{2}b^{2}}+\frac{3}{xb}+1\right)e^{-xb-ax^{2}}\, dx$?
The answer seems to be {}``no.'' Neither integration by parts, nor
positing some final form and checking that its derivative matches
the integrand yields a consistent result. Nor is Mathematica 13 able
to perform this indefinite integral. I checked that the only indefinite
integrals that Mathematica 13 is able to perform for a given \emph{$n\leq5$}
are those in Sections 2-5, apart from the known integrals having no
negative powers of \emph{x} in the integrand. Since the negative of
the indefinite integral with respect to \emph{a} of $b^{2}\int\frac{1}{x^{2}}\left(\frac{3}{x^{2}b^{2}}+\frac{3}{xb}+1\right)e^{-xb-ax^{2}}\, dx$,
gives $b^{2}\int\frac{1}{x^{4}}\left(\frac{3}{x^{2}b^{2}}+\frac{3}{xb}+1\right)e^{-xb-ax^{2}}\, dx$,
one might hope to find the negative of the indefinite integral with
respect to \emph{a} of the last line of (\ref{eq:V(-2,2,x) Theorem 6}),
but this, too, bore no fruit. Finally, one might hope to twice take
the indefinite integral with respect to $b$ of third and fourth lines
of (\ref{eq:V(-2,2,x) Theorem 6}), but the various factors of $b$
scattered throughout both the integrand and result foreclose this
option.

In a similar vein, since neither $b^{2}\int\left(\frac{3}{x^{2}b^{2}}+\frac{3}{xb}+1\right)e^{-xb-ax^{2}}\, dx$
nor $b^{2}\int\frac{1}{x^{2}}\left(\frac{3}{x^{2}b^{2}}+\frac{3}{xb}+1\right)e^{-xb-ax^{2}}\, dx$ 
poses a problem for any of the techniques we have used, one would
assume that the integrand with a power that lies between these two,
$b^{2}\int\frac{1}{x^{1}}\left(\frac{3}{x^{2}b^{2}}+\frac{3}{xb}+1\right)e^{-xb-ax^{2}}\, dx$,
would surely pose no difficulty since it has no larger pole at \emph{0}
than the former, and no higher value at infinity than the latter.
Unfortunately, that hope does not bear out for integration-by-hand
techniques introduced herein, nor can Mathematica do this integral.
The Rule-based Integration (Rubi) \cite{Rubi} package that augments
Mathematica likewise fails this task. The same negative result holds
for all other in-between integrals to those presented in Sections
2-5. The best one can do is a series solution for such integrals.

\subsection{The missing keystone integral}

One can find a series solution to any of these missing integrals from
the following theorem:

\subsection*{Theorem 16}

\subsection*{\begin{align}
\Upsilon\left(\upsilon,\, a,\, b,\, x\right) & =\int\frac{e^{-ax^{2}-bx}}{x^{\upsilon}}\, dx=\int\sqrt{\frac{2}{\pi}}e^{-ax^{2}}b^{0+\frac{1}{2}}x^{-0+\frac{1}{2}-\upsilon}K_{0+\frac{1}{2}}(xb)\, dx\nonumber \protect\\
= & \sum_{k=0}^{\infty}\left(\frac{b^{2k+1}a^{-k+\frac{\upsilon}{2}-1}\Gamma\left(k-\frac{\upsilon}{2}+1,\, ax^{2}\right)}{2(2k+1)!}-\frac{b^{2k}a^{-k+\frac{\upsilon}{2}-\frac{1}{2}}\Gamma\left(k-\frac{\upsilon}{2}+\frac{1}{2},\, ax^{2}\right)}{2(2k)!}\right)\:.\ \nonumber \protect\\
\equiv & \Upsilon_{e}\left(\upsilon,\, a,\, b,\, x\right)+\Upsilon_{o}\left(\upsilon,\, a,\, b,\, x\right)\label{eq:int exp-ax^2 series of cosh-sinh}\end{align}
 Proof. }

One may split $Exp(-bx)$ into even and odd parts, expand each term
separately,

\begin{align}
\Upsilon\left(\upsilon,\, a,\, b,\, x\right) & =\int\frac{e^{-ax^{2}}}{x^{\upsilon}}\left(\cosh(bx)-\sinh(bx)\right)\, dx\nonumber \\
= & \int\frac{e^{-ax^{2}}}{x^{\upsilon}}\left(\sum_{k=0}^{\infty}\frac{b^{2k}x^{2k}}{(2k)!}-\sum_{k=0}^{\infty}\frac{b^{2k+1}x^{2k+1}}{(2k+1)!}\right)\, dx\label{eq:cosh-sinh}\end{align}
 and integrate for any real or complex power $\upsilon$ of \emph{x}
to give (\ref{eq:int exp-ax^2 series of cosh-sinh}), thus completing
the proof. $\square$

The first five terms of the sum of this pair of series (evaluated
at $0.37$ minus the value at $0.31$) converge very rapidly: $\left.\Upsilon\left(1,\,0.11,\,0.13,x\right)\right|_{x\rightarrow0.31}^{x\rightarrow0.37}=\{0.1669998,\,0.000167671,\,2.76701\times10^{-8},\,1.83892\times10^{-12},\,6.60455\times10^{-17}\}$.
This is much more rapid than the alternative versions, below, and
avoids their alternating signs that might lead to roundoff errors.

The most obvious alternative approach is to expand %
{}

\begin{equation}
e^{-bx}=\sum_{k=0}^{\infty}\frac{(-bx)^{k}}{k!}\:\label{eq:exp-bx series}\end{equation}
 and integrate to obtain

\begin{equation}
\left.\Upsilon_{1}\left(1,\, a,\, b,\, x\right)\right|_{x\rightarrow c}^{x\rightarrow d}=\int_{c}^{d}\frac{e^{-ax^{2}-bx}}{x}\, dx=\int_{c}^{d}\frac{e^{-ax^{2}}}{x}\sum_{k=0}^{\infty}\frac{(-1)^{k}b^{k}x^{k}}{k!}\, dx=-\left.\sum_{k=0}^{\infty}\frac{(-1)^{k}a^{-\frac{k}{2}}b^{k}\Gamma\left(\frac{k}{2},\, ax^{2}\right)}{2k!}\right|_{x\rightarrow c}^{x\rightarrow d}\:.\label{eq:int exp-ax^2 series of exp-bx}\end{equation}

\noindent The first five terms of this series, again using the arbitrary
limits of integration $[c,d]=[0.31,0.37]$ are

\noindent $\left.\Upsilon_{1}\left(1,\,0.11,\,0.13,x\right)\right|_{x\rightarrow0.31}^{x\rightarrow0.37}\simeq\{0.174701,-0.0077012,\,0.000170185,-2.51374\times10^{-6},\,2.79188\times10^{-8}\}$,
and show decent but slower convergence than Theorem 16, and the alternating
signs may pose problems for other integral limit and parameter choices.
Or, as Nasser Abbasi suggested \cite{Abbasi}, one may instead expand
$e^{-ax^{2}}$ in a series, in which the above summand goes to $\frac{(-1)^{k}a^{k}b^{-2k}\Gamma(2k,bx)}{k!}$.
The first five terms for this version converge somewhat faster, $\left.\Upsilon_{2}\left(1,\,0.11,\,0.13,x\right)\right|_{x\rightarrow0.31}^{x\rightarrow0.37}\simeq\{0.169301,-0.00214673,\,0.0000137521,-5.933\times10^{-8},\,1.74623\times10^{-10}\}$,
but still have alternating signs.

There are dozens of  alternatives  to Theorem 16, particularly  insoluble one-dimensional
definite integrals, each ultimately requiring a series
solution (or numerical integration, which is a truncated series by
another name). The most notable of these (for $\upsilon=1$), involving
the incomplete plasma dispersion function, was introduced by Franklin
\cite{Franklin} and defined as

\begin{equation}
Z(u,w)\equiv\frac{1}{\sqrt{\pi}}\int_{u}^{\infty}\frac{e^{-t^{2}}}{t-w}dt\:.\label{eq:IPDF}\end{equation}

\noindent Baalrud~\cite{Baalrud} has provided MatLab code for computing
\emph{Z} via a Pade approximation or 10-point quadrature~\cite{Baalrud_code}.
Following the lead of Prudnikov, Brychkov, and Marichev \cite{PBM1}
(p. 140 No. 1.3.3.16) one changes variables to $x=\frac{2\sqrt{a}y-b}{2a}$
so that

\begin{equation}
\Upsilon\left(\upsilon,\, a,\, b,\, x\right)=e^{\frac{b^{2}}{4a}}a^{\frac{\nu}{2}-\frac{1}{2}}\int e^{-y^{2}}\left(y-\frac{b}{2\sqrt{a}}\right)^{-\nu}\, dy\:.\label{eq:exp(-x^2-x)/x-1}\end{equation}

\noindent Thus

\begin{eqnarray}
\left.\Upsilon_{Z}\left(1,\, a,\, b,\, x\right)\right|_{x\rightarrow c}^{x\rightarrow d}=\int_{c}^{d}\frac{e^{-ax^{2}-bx}}{x}\, dx=e^{\frac{b^{2}}{4a}}a^{\frac{\nu}{2}-\frac{1}{2}}\int_{\frac{2ac+b}{2\sqrt{a}}}^{\frac{2ad+b}{2\sqrt{a}}}e^{-y^{2}}\left(y-\frac{b}{2\sqrt{a}}\right)^{-1}\, dy\nonumber \\
=\sqrt{\pi}\left(Z\left(\frac{2ac+b}{2\sqrt{a}},\frac{b}{2\sqrt{a}}\right)-Z\left(\frac{2ad+b}{2\sqrt{a}},\frac{b}{2\sqrt{a}}\right)\right).\label{eq:int exp-ax^2 series of exp-bx-1}\end{eqnarray}

\noindent If the pole at $t=w$ in (\ref{eq:IPDF}) falls within the
range of integration ($u<w<\infty$), one subtracts an integral similar
to (\ref{eq:IPDF}) spanning $-\infty$ to $u$ from the full plasma
dispersion function \cite{fried,fadd}, $Z\left(-\infty,w\right)$.

There will be a numerical cost to using (\ref{eq:int exp-ax^2 series of cosh-sinh})
or any of the other four approximation methods, above, but none of
the other approaches I tried gave me an analytical result free of
infinte series.

In what follows, we will also need integrals with positive integer
powers of \emph{x} multiplying $e^{-ax^{2}-bx}$ that are given as
finite series in Prudnikov, Brychkov, and Marichev \cite{PBM1} (p.
140 No. 1.3.3.5,6,19). Their zero-power version is (p. 140 No. 1.3.3.17)
\begin{equation}
\Upsilon\left(0,a,b,x\right)=V\left(0,\,0,\, a,\, b,\, x\right)\equiv\int e^{-ax^{2}-bx}\, dx=\frac{\sqrt{\pi}e^{\frac{b^{2}}{4a}}\text{erf}\left(\frac{2ax+b}{2\sqrt{a}}\right)}{2\sqrt{a}}\:.\label{eq:unit-power}\end{equation}
 We will also want

\begin{equation}
V\left(1,\,0,\, a,\, b,\, x\right)\equiv\int xe^{-ax^{2}-bx}\, dx=-\frac{e^{-ax^{2}-bx}}{2a}-\frac{\sqrt{\pi}be^{\frac{b^{2}}{4a}}\text{erf}\left(\frac{2ax+b}{2\sqrt{a}}\right)}{4a^{3/2}}\:\label{eq:V[1,0]}\end{equation}
 and

\begin{equation}
V\left(2,\,0,\, a,\, b,\, x\right)\equiv\int x^{2}e^{-ax^{2}-bx}\, dx=\frac{e^{-ax^{2}-bx}(b-2ax)}{4a^{2}}+\frac{\sqrt{\pi}\left(2a+b^{2}\right)e^{\frac{b^{2}}{4a}}\text{erf}\left(\frac{2ax+b}{2\sqrt{a}}\right)}{8a^{5/2}}\:.\label{eq:V[2,0]}\end{equation}

\subsection{Usage}

Given that the integral of $e^{-ax^{2}}x^{-1+\frac{1}{2}}K_{\frac{3}{2}}(xb)$
can be done, shown in the last line of (\ref{eq:V(-1,1,x) Theorem 1}),
and one can likewise find the integral with all non-negative powers

\subsection*{Corollary 1}

\begin{equation}
\begin{array}{ll}
V\left(1,\,1,\, a,\, b,\, x\right)=\int\sqrt{\frac{2}{\pi}}e^{-ax^{2}}b^{\frac{3}{2}}x^{\frac{1}{2}+1}\frac{1}{2}G_{0,2}^{2,0}\left(\frac{b^{2}x^{2}}{4}|\begin{array}{c}
\frac{3}{4},-\frac{3}{4}\end{array}\right)\, dx\\
=\int\sqrt{\frac{2}{\pi}}e^{-ax^{2}}b^{\frac{3}{2}}x^{\frac{1}{2}+1}K_{\frac{3}{2}}(xb)\, dx\\
=b\int x\left(\frac{1}{bx}+1\right)e^{-ax^{2}-bx}\, dx=\frac{\sqrt{\pi}\left(2a-b^{2}\right)e^{\frac{b^{2}}{4a}}\text{erf}\left(\frac{2ax+b}{2\sqrt{a}}\right)}{4a^{3/2}}-\frac{be^{-ax^{2}-bx}}{2a}\end{array}\:\:\label{eq:V[1,1]}\end{equation}
 one would like to find the intermediate integral of $e^{-ax^{2}}x^{-0+\frac{1}{2}}K_{\frac{3}{2}}(xb)$.
No analytical result is to be found, but it can be given by the function
in the prior section:

\subsection*{Corollary 2}

\begin{equation}
\begin{array}{ll}
V\left(0,\,1,\, a,\, b,\, x\right)=\int\sqrt{\frac{2}{\pi}}\frac{b^{\frac{5}{2}}}{2^{2}}e^{-ax^{2}}x^{\frac{1}{2}+1}G_{0,2}^{2,0}\left(\frac{b^{2}x^{2}}{4}|\begin{array}{c}
\frac{1}{4},-\frac{5}{4}\end{array}\right)\, dx\\
=\int\sqrt{\frac{2}{\pi}}e^{-ax^{2}}b^{\frac{3}{2}}x^{\frac{1}{2}+0}K_{\frac{3}{2}}(xb)\, dx\\
=b\int\left(\frac{1}{bx}+1\right)e^{-ax^{2}-bx}\, dx=\Upsilon\left(1,a,b,x\right)+b\,\Upsilon\left(0,a,b,x\right)\end{array}\:\:\label{eq:V[0,1]}\end{equation}

\noindent From (\ref{eq:unit-power}) one sees that the final term
can also be written as $\Upsilon\left(0,a,b,x\right)\rightarrow V\left(0,\,0,\, a,\, b,\, x\right)$.

In between the integrands $e^{-ax^{2}}x^{-2+\frac{1}{2}}K_{\frac{5}{2}}(xb)$
and $e^{-ax^{2}}x^{0+\frac{1}{2}}K_{\frac{5}{2}}(xb)$ that pose no
problem for integration by parts, or checking derivatives of a posited
final form, or within Mathematica, one must likewise use the series
results to find the integral with a power intermediate to (\ref{eq:V(-2,2,x) Theorem 6})
and (\ref{eq:V(0,2,x) Theorem 2}):

\subsection*{Corollary 3}

\begin{equation}
\begin{array}{ll}
V\left(-1,\,2,\, a,\, b,\, x\right)=\sqrt{\frac{2}{\pi}}\int e^{-ax^{2}}b^{2+\frac{1}{2}}x^{-1+\frac{1}{2}}\frac{1}{2}G_{0,2}^{2,0}\left(\frac{b^{2}x^{2}}{4}|\begin{array}{c}
\frac{5}{4},-\frac{5}{4}\end{array}\right)\, dx\\
=\sqrt{\frac{2}{\pi}}\int e^{-ax^{2}}b^{2+\frac{1}{2}}x^{-1+\frac{1}{2}}K_{\frac{5}{2}}(xb)\,\, dx\\
=b^{2}\int\frac{1}{x}\left(\frac{3}{x^{2}b^{2}}+\frac{3}{xb}+1\right)e^{-xb-ax^{2}}\, dx\\
=3\Upsilon\left(3,a,b,x\right)+3b\Upsilon\left(2,a,b,x\right)+b^{2}\,\Upsilon\left(1,a,b,x\right)\:.\end{array}\:\label{eq:V(-1,2,x)}\end{equation}

\subsection{Higher-order integrals requiring series}

For $n=3$. the in-between the integrals are

\subsection*{Corollary 4}

\begin{equation}
\begin{array}{l}
V\left(0,\,3,\, a,\, b,\, x\right)=\sqrt{\frac{2}{\pi}}\int e^{-ax^{2}}b^{3+\frac{1}{2}}x^{0+\frac{1}{2}}\frac{1}{2}G_{0,2}^{2,0}\left(\frac{x^{2}b^{2}}{4}|\begin{array}{c}
\frac{7}{4},-\frac{7}{4}\end{array}\right)\, dx\\
=\sqrt{\frac{2}{\pi}}\int e^{-ax^{2}}b^{3+\frac{1}{2}}x^{0+\frac{1}{2}}K_{\frac{7}{2}}(xb)\, dx\\
=b^{3}\int\left(\frac{15}{x^{3}b^{3}}+\frac{15}{x^{2}b^{2}}+\frac{6}{xb}+1\right)e^{-xb-ax^{2}}\, dx\:\\
=15\Upsilon\left(3,a,b,x\right)+15b\Upsilon\left(2,a,b,x\right)+6b^{2}\,\Upsilon\left(1,a,b,x\right)+b^{3}V\left(0,\,0,\, a,\, b,\, x\right)\:\end{array}\:\label{eq:V(0,3,s)}\end{equation}

\noindent and

\subsection*{Corollary 5}

\begin{equation}
\begin{array}{l}
V\left(-2,\,3,\, a,\, b,\, x\right)=\int e^{-ax^{2}}b^{3+\frac{1}{2}}x^{-2+\frac{1}{2}}\frac{1}{2}G_{0,2}^{2,0}\left(\frac{x^{2}b^{2}}{4}|\begin{array}{c}
\frac{7}{4},-\frac{7}{4}\end{array}\right)\, dx\\
=\int e^{-ax^{2}}b^{-2+\frac{1}{2}}x^{-2+\frac{1}{2}}K_{\frac{7}{2}}(xb)\, dx\\
=b^{3}\int\frac{1}{x^{2}}\left(\frac{15}{x^{3}b^{3}}+\frac{15}{x^{2}b^{2}}+\frac{6}{xb}+1\right)e^{-ax^{2}-bx}\, dx\\
=15\Upsilon\left(5,a,b,x\right)+15b\Upsilon\left(4,a,b,x\right)+6b^{2}\,\Upsilon\left(3,a,b,x\right)+b^{3}\Upsilon\left(2,a,b,x\right)\end{array}\:.\label{eq:V(-2,3,s)}\end{equation}

There is no impediment to proceeding in this manner for any desired
value $n$ as long as one can keep track of terms that increase in
number as \emph{n} increases. Fortunately, a program like Mathematica,
can be programmed to access $V\left(h,\, n,\, a,\, b,\, x\right)$
and $\Upsilon\left(h,\, a,\, b,\, x\right)$ in an organized manner,
given in Appendix C. One finds that $n=4$ is sufficient to see the
trends in how the series is performing as the parameters vary.

For $n=4$, the in-between the integrals are

\subsection*{Corollary 6}

\begin{equation}
\begin{array}{ll}
V\left(1,\,4,\, a,\, b,\, x\right)=\sqrt{\frac{2}{\pi}}\int e^{-ax^{2}}b^{4+\frac{1}{2}}x^{1+\frac{1}{2}}\frac{1}{2}G_{0,2}^{2,0}\left(\frac{x^{2}b^{2}}{4}|\begin{array}{c}
\frac{9}{4},-\frac{9}{4}\end{array}\right)\, dx\\
=\sqrt{\frac{2}{\pi}}\int e^{-ax^{2}}b^{4+\frac{1}{2}}x^{1+\frac{1}{2}}K_{\frac{9}{2}}(xb)\,\, dx\:\\
=b^{4}\int x\left(\frac{105}{x^{4}b^{4}}+\frac{105}{x^{3}b^{3}}+\frac{45}{x^{2}b^{2}}+\frac{10}{xb}+1\right)e^{-ax^{2}-bx}\,\, dx\\
=105\Upsilon\left(3,a,b,x\right)+105b\Upsilon\left(2,a,b,x\right)+45b^{2}\,\Upsilon\left(1,a,b,x\right)\\
+10b^{3}\,\Upsilon\left(0,a,b,x\right)\:+b^{4}V\left(1,\,0,\, a,\, b,\, x\right)\end{array}\:.,\label{eq:V(1,4,s)}\end{equation}
(where from (\ref{eq:unit-power}) the penultimate term can also be
written as $10b^{3}\,\Upsilon\left(0,a,b,x\right)\rightarrow10b^{3}\, V\left(0,\,0,\, a,\, b,\, x\right)$),

\subsection*{Corollary 7}

\begin{equation}
\begin{array}{ll}
V\left(-1,\,4,\, a,\, b,\, x\right)=\sqrt{\frac{2}{\pi}}\int e^{-ax^{2}}b^{4+\frac{1}{2}}x^{-1+\frac{1}{2}}\frac{1}{2}G_{0,2}^{2,0}\left(\frac{x^{2}b^{2}}{4}|\begin{array}{c}
\frac{9}{4},-\frac{9}{4}\end{array}\right)\, dx\\
=\sqrt{\frac{2}{\pi}}\int e^{-ax^{2}}b^{4+\frac{1}{2}}x^{-1+\frac{1}{2}}K_{\frac{9}{2}}(xb)\,\, dx\:\\
=b^{4}\int\frac{1}{x}\left(\frac{105}{x^{4}b^{4}}+\frac{105}{x^{3}b^{3}}+\frac{45}{x^{2}b^{2}}+\frac{10}{xb}+1\right)e^{-ax^{2}-bx}\,\, dx\\
=105\Upsilon\left(5,a,b,x\right)+105b\Upsilon\left(4,a,b,x\right)+45b^{2}\,\Upsilon\left(3,a,b,x\right)\:\\
+10b^{3}\Upsilon\left(2,a,b,x\right)+b^{4}\Upsilon\left(1,a,b,x\right)\end{array}\:,\label{eq:V(-1,4,s)}\end{equation}
and

\subsection*{Corollary 8}

\begin{equation}
\begin{array}{ll}
V\left(-3,\,4,\, a,\, b,\, x\right)=\sqrt{\frac{2}{\pi}}\int e^{-ax^{2}}b^{4+\frac{1}{2}}x^{-3+\frac{1}{2}}\frac{1}{2}G_{0,2}^{2,0}\left(\frac{x^{2}b^{2}}{4}|\begin{array}{c}
\frac{9}{4},-\frac{9}{4}\end{array}\right)\, dx\\
=\sqrt{\frac{2}{\pi}}\int e^{-ax^{2}}b^{4+\frac{1}{2}}x^{-3+\frac{1}{2}}K_{\frac{9}{2}}(xb)\,\, dx\:\\
=b^{4}\int\frac{1}{x^{3}}\left(\frac{105}{x^{4}b^{4}}+\frac{105}{x^{3}b^{3}}+\frac{45}{x^{2}b^{2}}+\frac{10}{xb}+1\right)e^{-ax^{2}-bx}\,\, dx\\
=105\Upsilon\left(7,a,b,x\right)+105b\Upsilon\left(6,a,b,x\right)+45b^{2}\,\Upsilon\left(5,a,b,x\right)\:\\
+10b^{3}\Upsilon\left(4,a,b,x\right)+b^{4}\Upsilon\left(3,a,b,x\right)\end{array}\:.\label{eq:V(-3,4,s)}\end{equation}

One can follow the pattern established for $n=4$ to write down the
in-between the integrals for $n=5$.

\section{Incomplete Gamma Function}

Tables of incomplete gamma functions are even more spare than for
Bessel functions. Indeed Prudnikov, Brychkov, and Marichev~\cite{PBM2}
(p. 23 Section 1.2.2 containiing the lower incomplete gamma function
$\gamma(j,bx)\equiv\Gamma(j)-\Gamma(j,bx)$) contains precisely four
involving exponentials and just one of those is multiplied by any
power power ($x^{1}$). Of upper incomplete gamma functions with integer
indices, the only integral of the form

\begin{equation}
\int\frac{e^{-ax^{2}}\Gamma(j,bx)}{x^{\nu}}\, dx\label{eq:only_Gamma_int}\end{equation}

\noindent that Mathematica 13 can do is when $\nu=0$ and $j=1$,
which integrand is the pure exponential $e^{-ax^{2}-bx}$, unless
-- also for $\nu=0$ -- one expands $\Gamma(j,bx)$ in $e^{-bx}$
multiplied by positive powers~\cite{GR5} (p. 949 No. 8.352.2).

Integrals for $\nu>0$ may be found by inserting the set of theorems
so far presented into the following:

\subsection*{Theorem 17}

\subsection*{\begin{eqnarray}
&& \hspace{-1.1 cm}  \int\frac{ e^{-ax^{2}} \Gamma(2L+1+m,bx)}{b^{2L+1}x^{2L+1}}\, dx = \nonumber \protect\\
& = & (2L+m)!\sqrt{\frac{2}{\pi}}\;\left(\frac{1}{(2L-1)!!}\;\sum_{j=1}^{L}\;\frac{1}{[2(L-j)]!!}\int  e^{-ax^{2}}  b^{-(j+1/2)}x^{-(j+1/2)}\; K_{j+1/2}(bx)\, dx\right.\nonumber \protect\\
 & + & \left.\frac{1}{(2L)!}\int  e^{-ax^{2}}  b^{-(1/2)}x^{-(1/2)}\; K_{1/2}(bx)\, dx\, dx\right)\;,\label{eq: int of twentytwo}\protect\\
 & + & \sum_{h=0}^{m-1}\;\sum_{j=h}^{m-1}\;\frac{a(j,\, h)}{(2L+1+j)!}\int  e^{-ax^{2}}  b^{h+1/2}x^{h+1/2}\; K_{h+1/2}(bx) \nonumber 
 \end{eqnarray}
Proof. }

In 2001~\cite{stra01} I showed that the incomplete gamma function
can be expanded in a finite series of Macdonald functions,

\begin{eqnarray}
\frac{\Gamma(2L+1+m,\zeta)}{\zeta^{2L+1}} & = & (2L+m)!\sqrt{\frac{2}{\pi}}\;\left(\frac{1}{(2L-1)!!}\;\sum_{j=1}^{L}\;\frac{\zeta^{-(j+1/2)}\; K_{j+1/2}(\zeta)}{[2(L-j)]!!}\right.\nonumber \\
 & + & \left.\frac{\zeta^{-(1/2)}\; K_{1/2}(\zeta)}{(2L)!}+\sum_{h=0}^{m-1}\;\sum_{j=h}^{m-1}\;\frac{a(j,\, h)\zeta^{h+1/2}\; K_{h+1/2}(\zeta)}{(2L+1+j)!}\right)\;,\label{twentytwo}\end{eqnarray}

\noindent where the $a$'s are given by recursion, \begin{eqnarray}
a(i,\, i) & \equiv & 1\nonumber \\
a(i,\, i-1) & = & -\left({i+1\atop 2}\right)\nonumber \\
a(i,\, i-k) & = & -\sum_{m=0}^{k-1}a_{i}^{i-m}\left({i+k-2m\atop 2k-2m}\right)(2k-2m-1)!!\nonumber \\
a(i,\, k) & \equiv & 0\;\;,k<\left\lfloor (i-1)/2\right\rfloor \;.\label{fourteen}\end{eqnarray}

\noindent Unless $k$ is greater than 1 and less than the greatest
integer less than or equal to $(i+1)/2$, symbolically written $1<k\leq\left\lfloor \left(i+1\right)/2\right\rfloor $,
the third line gives zero-valued results. To save work, one should
just skip to the fourth line without calculating it, though zero results
in either case. In that prior paper~\cite{stra01} I had put such
a redundant advisory note on the third line that erroneously said
that unless $1<k<[i/2]$, one should just skip to the fourth line.
The first eight sets of $a$'s are\begin{equation}
\begin{array}{c}
\{a(0,\,0)\to1\}\\
\{a(1,\,1)\to1,\, a(1,\,0)\to-1\}\\
\{a(2,\,2)\to1,\, a(2,\,1)\to-3,\, a(2,\,0)\to0\}\\
\{a(3,\,3)\to1,\, a(3,\,2)\to-6,\, a(3,\,1)\to3,\, a(3,\,0)\rightarrow0\}\\
\{a(4,\,4)\to1,\, a(4,\,3)\to-10,\, a(4,\,2)\to15,\, a(4,\,1)=a(4,\,0)\to0\}\\
\{a(5,\,5)\to1,\, a(5,\,4)\to-15,\, a(5,\,3)\to45,\, a(5,\,2)\to-15.a(5,\,0)\to0=a(5,\,1)\to0\}\\
\{a(6,\,6)\to1,\, a(6,\,5)\to-21,\, a(6,\,4)\to105,\, a(6,\,3)\to-105,\, a(6,\,2)\to0,\, a(6,\,1)=a(6,\,0)\to0\}\\
\text{\footnotesize\ensuremath{\{a(7,\,7)\to1,\, a(7,\,6)\to-28,\, a(7,\,5)\to210,\, a(7,\,4)\to-420,\, a(7,\,3)\to105,\, a(7,\,2)=a(7,\,1)=a(7,\,0)\to0\}}}\\
\text{\footnotesize\ensuremath{\{a(8,\,8)\to1,\, a(8,\,7)\to-36,\, a(8,\,6)\to378,\, a(8,\,5)\to-1260,\, a(8,\,4)\to945,\, a(8,\,3)=a(8,\,2)=a(8,\,1)=a(8,\,0)\to0\}}}\end{array}\label{eq:a[j,h]}\end{equation}

\noindent One may integrate (\ref{twentytwo}) term-by-term to complete
the proof. $\square$ 

This allows one to find the indefinite integral of many incomplete
Gamma functions combined with inverse powers and $e^{-ax^{2}}$. For
$L=m=0$ the integral is

\subsection*{Corollary 9\begin{equation}
\int\frac{e^{-ax^{2}}\Gamma(1,bx)}{x}\, dx=\int\frac{e^{-ax^{2}-bx}}{x}\, dx=\Upsilon(1,\, a,\, b,\, x)\:.\label{eq:Int_Gamma1/x}\end{equation}
}

\noindent A numerical check of the above integral over the interval
$\left[0.31,0.37\right]$, gives the result $0.16716752839055965$,
and five terms in the series for $\Upsilon(1,\,0.11,\,0.13,\, x)$
are sufficent for seventeen digit accuracy.

For $L=0$ and $m=1$ the integral is a mix of \emph{V}'s and $\Upsilon$'s:

\subsection*{Corollary 10}

\begin{equation}
\begin{array}{ccc}
\int\frac{e^{-ax^{2}}\Gamma(2,bx)}{x}\, dx & = & \int b\left(\frac{\sqrt{\frac{2}{\pi}}e^{-ax^{2}}K_{\frac{1}{2}}(bx)}{\sqrt{bx}}+0+\left\{ \sqrt{\frac{2}{\pi}}a(0,0)e^{-ax^{2}}\sqrt{bx}K_{\frac{1}{2}}(bx)\right\} \right)\, dx\\
 & = & \int\left(\frac{e^{-ax^{2}-bx}}{x}+be^{-ax^{2}-bx}\right)\, dx\\
 & = & \Upsilon(1,\, a,\, b,\, x)+bV(0,\,0,\, a,\, b,\, x)\end{array}\label{eq:int_Gamma2/x}\end{equation}
 For higher values of \emph{m}, such as $L=0$ and $m=2$, one also
has a mix of $\Upsilon$'s and \emph{V}'s, the latter being integrals
of positive integer powers of \emph{x} multiplying $e^{-ax^{2}-bx}$
given in Prudnikov, Brychkov, and Marichev, \cite{PBM1} (p. 140 No.
1.3.3.5,6,19):

\subsection*{Corollary 11}

\noindent \begin{equation}
\begin{array}{ccc}
\int\frac{e^{-ax^{2}}\Gamma(3,bx)}{x}\, dx & = & \int b\left(\frac{2\sqrt{\frac{2}{\pi}}e^{-ax^{2}}K_{\frac{1}{2}}(bx)}{\sqrt{bx}}+0+\sqrt{\frac{2}{\pi}}a(1,1)e^{-ax^{2}}(bx)^{3/2}K_{\frac{3}{2}}(bx)\right.\\
 & + & \left.2\sqrt{\frac{2}{\pi}}e^{-ax^{2}}\left(a(0,0)\sqrt{bx}K_{\frac{1}{2}}(bx)+\frac{1}{2}a(1,0)\sqrt{bx}K_{\frac{1}{2}}(bx)\right)\right)\, dx\\
 & = & \int b\left(2\frac{e^{-ax^{2}-bx}}{x}+e^{-ax^{2}-bx}+bx\left(\frac{1}{bx}+1\right)e^{-ax^{2}-bx}\right)\, dx\\
 & = & 2\Upsilon(1,\, a,\, b,\, x)+bV(0,\,0,\, a,\, b,\, x)+bV(1,\,1,\, a,\, b,\, x)\end{array}\label{eq:int_Gamma3/x}\end{equation}
where the final term in line three is $b\times b\int x\left(\frac{1}{bx}+1\right)e^{-ax^{2}-bx}\, dx=b\times V\left(1,\,1,\, a,\, b,\, x\right)$.
In moving from line three to four, one must remember that each $V(h,\, n,\, a,\, b,\, x)$
already contains a factor $b^{n}$ but the definition of $\Upsilon(\upsilon,\, a,\, b,\, x)$
does not.

In all integrals that follow, instead of displaying the various $a(i,k)$,
like those in the second line, above, we will insert their values
and sum common terms. For $L=1$ and $m=0$ 

\subsection*{Corollary 12}

\begin{equation}
\begin{array}{ccc}
\int\frac{e^{-ax^{2}}\Gamma(3,bx)}{b^{2}x^{3}}\, dx & = & \int\left(\frac{\sqrt{\frac{2}{\pi}}be^{-ax^{2}}K_{\frac{1}{2}}(bx)}{\sqrt{bx}}+\frac{2\sqrt{\frac{2}{\pi}}be^{-ax^{2}}K_{\frac{3}{2}}(bx)}{(bx)^{3/2}}+0\right)\, dx\\
 & = & \int\left(\frac{e^{-ax^{2}-bx}}{x}+\frac{2\left(1+\frac{1}{bx}\right)e^{-ax^{2}-bx}}{bx^{2}}\right)\, dx\:.\\
 & = & \Upsilon(1,\, a,\, b,\, x)+\frac{2}{b}\Upsilon(2,\, a,\, b,\, x)+\frac{2}{b^{2}}\Upsilon(3,\, a,\, b,\, x)\end{array}\label{eq:int_Gamma3/x^3}\end{equation}
 With $m=1$,

\subsection*{Corollary 13}

\noindent \begin{equation}
\begin{array}{ccc}
\int\frac{e^{-ax^{2}}\Gamma(4,bx)}{b^{2}x^{3}}\, dx & = & \int\left(\frac{3\sqrt{\frac{2}{\pi}}be^{-ax^{2}}K_{\frac{1}{2}}(bx)}{\sqrt{bx}}+\frac{6\sqrt{\frac{2}{\pi}}be^{-ax^{2}}K_{\frac{3}{2}}(bx)}{(bx)^{3/2}}\right.\\
 & + & \left.\sqrt{\frac{2}{\pi}}be^{-ax^{2}}\sqrt{bx}K_{\frac{1}{2}}(bx)\right)\, dx\:,\\
 & = & \int\left(\frac{3e^{-ax^{2}-bx}}{x}+\frac{6\left(1+\frac{1}{bx}\right)e^{-ax^{2}-bx}}{bx^{2}}+be^{-ax^{2}-bx}\right)\, dx\\
 & = & 3\Upsilon(1,\, a,\, b,\, x)+\frac{6}{b}\Upsilon(2,\, a,\, b,\, x)+\frac{6}{b^{2}}\Upsilon(3,\, a,\, b,\, x)+bV(0,\,0,\, a,\, b,\, x)\end{array}\label{eq:int_Gamma4/x^3}\end{equation}
 and with $m=2$,

\subsection*{Corollary 14}

\noindent \begin{equation}
\begin{array}{ccc}
 \hspace{-0.3 cm} \int\frac{e^{-ax^{2}}\Gamma(5,bx)}{b^{2}x^{3}}\,  dx \hspace{-0.2 cm} & = &\hspace{-0.1cm}    \int\left(\frac{12\sqrt{\frac{2}{\pi}}be^{-ax^{2}}K_{\frac{1}{2}}(bx)}{\sqrt{bx}}+\frac{24\sqrt{\frac{2}{\pi}}be^{-ax^{2}}K_{\frac{3}{2}}(bx)}{(bx)^{3/2}}\right.\\
 & + & \left.3\sqrt{\frac{2}{\pi}}be^{-ax^{2}}\sqrt{bx}K_{\frac{1}{2}}(bx)+\sqrt{\frac{2}{\pi}}be^{-ax^{2}}(bx)^{3/2}K_{\frac{3}{2}}(bx)\right)\, dx\:.\\
 & = & \int\left(\frac{12e^{-ax^{2}-bx}}{x}+\frac{24\left(\frac{1}{bx}+1\right)e^{-ax^{2}-bx}}{bx^{2}}+3be^{-ax^{2}-bx}+b^{2}x\left(\frac{1}{bx}+1\right)e^{-ax^{2}-bx}\right)\, dx\\
 & = & \hspace{-0.3cm}12\Upsilon(1,\, a,\, b,\, x)+\frac{24}{b}\Upsilon(2,\, a,\, b,\, x)+\frac{24}{b^{2}}\Upsilon(3,\, a,\, b,\, x)+3bV(0,\,0,\, a,\, b,\, x)+bV(1,\,1,\, a,\, b,\, x)\end{array}\label{eq:int_Gamma5/x^3}\end{equation}
Five terms in the series for $\Upsilon(n,\,0.11,\,0.13,\, x)$ are
sufficient for sixteen digit accuracy in comparison with numerical
integration checks of the above three integrals over the interval
$\left[0.31,0.37\right]$.

Moving to $L=2$ and $m=0$

\subsection*{Corollary 15}

\begin{equation}
\begin{array}{ccc}
\int\frac{e^{-ax^{2}}\Gamma(5,bx)}{b^{4}x^{5}}\, dx & = & \int\left(\frac{\sqrt{\frac{2}{\pi}}be^{-ax^{2}}K_{\frac{1}{2}}(bx)}{\sqrt{bx}}+\frac{4\sqrt{\frac{2}{\pi}}be^{-ax^{2}}K_{\frac{3}{2}}(bx)}{(bx)^{3/2}}+\frac{8\sqrt{\frac{2}{\pi}}be^{-ax^{2}}K_{\frac{5}{2}}(bx)}{(bx)^{5/2}}\right)\, dx\\
 & = & \int\left(\frac{e^{-ax^{2}-bx}}{x}+\frac{4\left(1+\frac{1}{bx}\right)e^{-ax^{2}-bx}}{bx^{2}}+\frac{8\left(1+\frac{3}{bx}+\frac{3}{b^{2}x^{2}}\right)e^{-ax^{2}-bx}}{b^{2}x^{3}}\right)\, dx\:.\\
 & = & \Upsilon(1,\, a,\, b,\, x)+\frac{4}{b}\Upsilon(2,\, a,\, b,\, x)+\frac{4+8}{b^{2}}\Upsilon(3,\, a,\, b,\, x)\\
 & + & \frac{24}{b^{3}}\Upsilon(4,\, a,\, b,\, x)+\frac{24}{b^{4}}\Upsilon(5,\, a,\, b,\, x)\end{array}\label{eq:int_Gamma5/x^5}\end{equation}

\noindent With $m=1$,

\subsection*{Corollary 16}

\begin{equation}
\begin{array}{ccc}
\int\frac{e^{-ax^{2}}\Gamma(6,bx)}{b^{4}x^{5}}\, dx & = & \int\left(\frac{5\sqrt{\frac{2}{\pi}}be^{-ax^{2}}K_{\frac{1}{2}}(bx)}{\sqrt{bx}}+\frac{20\sqrt{\frac{2}{\pi}}be^{-ax^{2}}K_{\frac{3}{2}}(bx)}{(bx)^{3/2}}+\frac{40\sqrt{\frac{2}{\pi}}be^{-ax^{2}}K_{\frac{5}{2}}(bx)}{(bx)^{5/2}}\right.\\
 & + & \left.\sqrt{\frac{2}{\pi}}be^{-ax^{2}}\sqrt{bx}K_{\frac{1}{2}}(bx)\right)\, dx\:,\\
 & = & \int\left(\frac{5e^{-ax^{2}-bx}}{x}+\frac{20\left(\frac{1}{bx}+1\right)e^{-ax^{2}-bx}}{bx^{2}}+\frac{40\left(\frac{3}{b^{2}x^{2}}+\frac{3}{bx}+1\right)e^{-ax^{2}-bx}}{b^{2}x^{3}}+be^{-ax^{2}-bx}\right)\, dx\\
 & = & \hspace{-1.1cm}5\Upsilon(1,a,b,x)+\frac{20}{b}\Upsilon(2,a,b,x)+\frac{60}{b^{2}}\Upsilon(3,a,b,x)\\
 & + & \frac{120}{b^{3}}\Upsilon(4,a,b,x)+\frac{120}{b^{4}}\Upsilon(5,a,b,x)+bV(0,0,a,b,x)\end{array}\label{eq:int_Gamma6/x^5}\end{equation}
 and with $m=2$,

\subsection*{Corollary 17}

\noindent \begin{equation}
\begin{array}{ccc}
\int\frac{e^{-ax^{2}}\Gamma(7,bx)}{b^{4}x^{5}}\, dx\hspace{-0.3cm} & = & \int\left(\frac{30\sqrt{\frac{2}{\pi}}be^{-ax^{2}}K_{\frac{1}{2}}(bx)}{\sqrt{bx}}+\frac{120\sqrt{\frac{2}{\pi}}be^{-ax^{2}}K_{\frac{3}{2}}(bx)}{(bx)^{3/2}}+\frac{240\sqrt{\frac{2}{\pi}}be^{-ax^{2}}K_{\frac{5}{2}}(bx)}{(bx)^{5/2}}\right.\\
 & + & \left.5\sqrt{\frac{2}{\pi}}be^{-ax^{2}}\sqrt{bx}K_{\frac{1}{2}}(bx)+\sqrt{\frac{2}{\pi}}be^{-ax^{2}}(bx)^{3/2}K_{\frac{3}{2}}(bx)\right)\, dx\:.\\
 & = & \int\left(\frac{30e^{-ax^{2}-bx}}{x}+\frac{120\left(\frac{1}{bx}+1\right)e^{-ax^{2}-bx}}{bx^{2}}+\frac{240\left(\frac{3}{b^{2}x^{2}}+\frac{3}{bx}+1\right)e^{-ax^{2}-bx}}{b^{2}x^{3}}\right.\\
 & + & \left.5be^{-ax^{2}-bx}+b^{2}x\left(\frac{1}{bx}+1\right)e^{-ax^{2}-bx}\right)\, dx\\
 & = & 30\Upsilon(1,\, a,\, b,\, x)+\frac{120}{b}\Upsilon(2,\, a,\, b,\, x)+\frac{360}{b^{2}}\Upsilon(3,\, a,\, b,\, x)\\
 & + & \hspace{-0.3cm}\frac{720}{b^{3}}\Upsilon(4,\, a,\, b,\, x)+\frac{720}{b^{4}}\Upsilon(5,\, a,\, b,\, x)+5bV(0,\,0,\, a,\, b,\, x)+bV(1,\,1,\, a,\, b,\, x)\end{array}\label{eq:int_Gamma7/x^5}\end{equation}
For the above three integrals, five terms in the series for $\Upsilon(n,\,0.11,\,0.13,\, x)$
are sufficient to match numerical integration over the interval $\left[0.31,0.37\right]$
to eighteen digit accuracy.

\subsection{Integrals with only analytical functions in the result}

\noindent The relation (\ref{twentytwo}) between incomplete gamma
functions and Macdonald functions of course holds if one multiplies
both sides by a power of $\zeta=bx$, so this can be used to find
indefinite integrals with alternative powers for a given gamma function.
For instance, if we multiply the integrand of (\ref{eq:int_Gamma3/x^3})
by $bx$, the resulting integral \textendash{} still with $L=1$ and
$m=0$ \textendash{} has only analytical functions in the result:

\subsection*{Corollary 18}

\begin{equation}
\begin{array}{ccc}
\int\frac{e^{-ax^{2}}\Gamma(3,bx)}{bx^{2}}\, dx & = & \int\left(be^{-ax^{2}-bx}+\frac{2\left(\frac{1}{bx}+1\right)e^{-ax^{2}-bx}}{x}\right)\, dx\\
 & = & bV(0,\,0,\, a,\, b,\, x)+\frac{2}{b}V(-1,\,1,\, a,\, b,\, x)\:.\end{array}\label{eq:int_Gamma3/x^2_is_anyl}\end{equation}
 Similarly, multiplying the integrand of (\ref{eq:int_Gamma4/x^3}),
derived from (\ref{twentytwo}) with $L=1$ and $m=1$ , gives

\subsection*{Corollary 19}

\noindent \begin{equation}
\begin{array}{ccc}
\int\frac{e^{-ax^{2}}\Gamma(4,bx)}{bx^{2}}\, dx & = & \int\left(3be^{-ax^{2}-bx}+\frac{6\left(\frac{1}{bx}+1\right)e^{-ax^{2}-bx}}{x}+b^{2}xe^{-ax^{2}-bx}\right)\, dx\\
 & = & 3bV(0,\,0,\, a,\, b,\, x)+\frac{6}{b}V(-1,\,1,\, a,\, b,\, x)+b^{2}V(1,\,0,\, a,\, b,\, x)\:.\end{array}\label{eq:int_Gamma4/x^2_is_anyl}\end{equation}
The $L=1$ and $m=2$ version will be found in the following section.

If we multiply (\ref{eq:int_Gamma5/x^5}) by $bx$, the resulting
integral \textendash{} still with $L=2$ and $m=0$ \textendash{}
likewise has only analytical functions in the result,

\subsection*{Corollary 20}

\begin{equation}
\begin{array}{ccc}
\int\frac{e^{-ax^{2}}\Gamma(5,bx)}{b^{3}x^{4}}\, dx & = & \int\left(be^{-ax^{2}-bx}+\frac{4\left(\frac{1}{bx}+1\right)e^{-ax^{2}-bx}}{x}+\frac{8\left(\frac{3}{b^{2}x^{2}}+\frac{3}{bx}+1\right)e^{-ax^{2}-bx}}{bx^{2}}\right)\, dx\\
 & = & bV(0,\,0,\, a,\, b,\, x)+\frac{4}{b}V(-1,\,1,\, a,\, b,\, x)+\frac{8}{b^{3}}V(-2,\,2,\, a,\, b,\, x)\:,\end{array}\label{eq:int_Gamma5/x^4}\end{equation}
 as does (\ref{eq:int_Gamma6/x^5}) (which has $L=2$ and $m=1$)
if we multiply by $bx$:

\subsection*{Corollary 21}

\noindent \begin{equation}
\begin{array}{ccc}
\int\frac{e^{-ax^{2}}\Gamma(6,bx)}{b^{3}x^{4}}\, dx & = & \int\left(5be^{-ax^{2}-bx}+\frac{20\left(\frac{1}{bx}+1\right)e^{-ax^{2}-bx}}{x}+\frac{40\left(\frac{3}{b^{2}x^{2}}+\frac{3}{bx}+1\right)e^{-ax^{2}-bx}}{bx^{2}}+b^{2}xe^{-ax^{2}-bx}\right)\, dx\\
 & = & 5bV(0,\,0,\, a,\, b,\, x)+\frac{20}{b}V(-1,\,1,\, a,\, b,\, x)+\frac{40}{b^{3}}V(-2,\,2,\, a,\, b,\, x)+b^{2}V(1,\,1,\, a,\, b,\, x)\:.\end{array}\label{eq:int_Gamma6/x^4}\end{equation}
 Analytical results are also obtained if we multiply the integrand
of (\ref{eq:int_Gamma7/x^5}) (which has $L=2$ and $m=2$) by $bx$:

\subsection*{Corollary 22}

\noindent \begin{equation}
\begin{array}{ccc}
\int\frac{e^{-ax^{2}}\Gamma(7,bx)}{b^{3}x^{4}}\, dx & = & \int\left(30be^{-ax^{2}-bx}+\frac{120\left(b\left(\left(\frac{1}{bx}+1\right)e^{-ax^{2}-bx}\right)\right)}{bx}+\frac{240\left(b^{2}\left(\left(\frac{3}{b^{2}x^{2}}+\frac{3}{bx}+1\right)e^{-ax^{2}-bx}\right)\right)}{b^{3}x^{2}}\right.\\
 & + & \left.5b^{2}\left(xe^{-ax^{2}-bx}\right)+b^{2}\left(bx^{2}\left(\frac{1}{bx}+1\right)e^{-ax^{2}-bx}\right)\right)\, dx \;.\\
 & = & 30bV(0,\,0,\, a,\, b,\, x)+\frac{120}{b}V(-1,\,1,\, a,\, b,\, x)+\frac{240}{b^{3}}V(-2,\,2,\, a,\, b,\, x)\\
 & + & 5b^{2}V(1,\,0,\, a,\, b,\, x)\:+bV(2,1,\, a,\, b,\, x)\end{array}\label{eq:int_Gamma7/x^4}\end{equation}

It is straightforward to continue to use (\ref{eq: int of twentytwo})
to obtain integrals with larger values of \emph{j} in $\Gamma(j,\, bx)$,
for both series and analytical functions in the result, but an alternative
route is also open.

\subsection{Recursion approach}

One may apply the recursion relation \cite{GR5} (p. 951 No. 8.356.2)
\begin{equation}
\Gamma(a+1,z)=e^{-z}z^{a}+a\Gamma(a,z)\label{eq:recursion}\end{equation}
in the simplest case by using the the known integrals (\ref{eq:Int_Gamma1/x})
and (\ref{eq:unit-power}) after stepping down the left-hand side
of (\ref{eq:int_Gamma2/x}) to again obtain

\subsection*{Corollary 23}

\noindent \begin{equation}
\begin{array}{ccc}
\int\frac{e^{-ax^{2}}\Gamma(2,bx)}{x}\, dx & = & \int\frac{e^{-ax^{2}}}{x}\left(\Gamma(1,bx)+bxe^{-bx}\right)\, dx=\Upsilon(1,\, a,\, b,\, x)+\int\frac{e^{-ax^{2}}}{x}\left(bxe^{-bx}\right)\, dx\\
 & = & \Upsilon(1,\, a,\, b,\, x)+bV(0,\,0,\, a,\, b,\, x)\end{array}\label{eq:recursionGamma2/x}\end{equation}

Of more interest is to obtain integrals not listed in the above sections,
particularly those with analytical solutions rather than series solutions.
One may re-arrange (\ref{eq:recursion}) to give a step-up version,
\begin{equation}
\Gamma(a-1,\, z)=\frac{1}{a-1}\left(\Gamma(a,\, z)-e^{-z}z^{a-1}\right)\:.\label{eq:step-up}\end{equation}
 Let us apply this to the case where we multiplied the $L=1$ and
$m=0$ integrand (\ref{eq:int_Gamma3/x^3}) by $bx$, which gave us
(\ref{eq:int_Gamma3/x^2_is_anyl}). One can step up from the integrand
on the left-hand side of the following, yielding the integral:

\subsection*{Corollary 24}

\noindent \begin{equation}
\begin{array}{ccc}
\int\frac{e^{-ax^{2}}\Gamma(2,bx)}{bx^{2}}\, dx & = & \int\frac{e^{-ax^{2}}}{bx^{2}}\left(\frac{1}{2}\Gamma(3,bx)-\frac{1}{2}b^{2}x^{2}e^{-bx}\right)\, dx\\
 & = & \frac{1}{b}V(-1,\,1,\, a,\, b,\, x)\:.\end{array}\label{eq:int_Gamma2/x^2_is_anyl}\end{equation}
where the $-\frac{1}{2}bV(0,\,0,\, a,\, b,\, x)$ from the second
term 
on the right-hand side of the first line cancels the $\frac{1}{2}bV(0,\,0,\, a,\, b,\, x)$
from the first term in (\ref{eq:int_Gamma3/x^2_is_anyl}). Again this
integral has only analytical functions in the result.

Stepping up again provides an alternative way to obtain (\ref{eq:int_Gamma4/x^2_is_anyl})
from the prior section, and two more steps gives

\subsection*{Corollary 25}

\begin{equation}
\begin{array}{ccc}
\int\frac{e^{-ax^{2}}\Gamma(5,bx)}{bx^{2}}\, dx & = & 12bV(0,\,0,\, a,\, b,\, x)+\frac{24}{b}V(-1,\,1,\, a,\, b,\, x)+4b^{2}V(1,\,0,\, a,\, b,\, x)+b^{3}V(2,\,0,\, a,\, b,\, x)\:.\end{array}\label{eq:int_Gamma5/x^2_is_anyl}\end{equation}

Staying with the case where $L=1$ and $m=0$ but multiplying the
integrand (\ref{eq:int_Gamma3/x^3}) by $b^{3}x^{3}$ again gives
analytical functions in the result,

\subsection*{Corollary 26}

\begin{equation}
\begin{array}{ccc}
\int e^{-ax^{2}}\Gamma(3,bx)\, dx & = & b^{2}V(2,\,0,\, a,\, b,\, x)+2V(1,\,1,\, a,\, b,\, x)\end{array}\label{eq:int_Gamma3}\end{equation}
 and stepping down with (\ref{eq:step-up}) gives,

\subsection*{Corollary 27}

\noindent \begin{equation}
\begin{array}{ccc}
\int e^{-ax^{2}}\Gamma(2,bx)\, dx & = & V(1,\,1,\, a,\, b,\, x)\;.\end{array}\label{eq:int_Gamma2}\end{equation}

For the case where $L=2$ and $m=0$ but multiplying the integrand
(\ref{eq:int_Gamma5/x^5}) by $b^{3}x^{3}$ gives

\subsection*{Corollary 28}

\begin{equation}
\begin{array}{ccc}
\int\frac{e^{-ax^{2}}\Gamma(5,bx)}{bx^{2}}\, dx & =b^{3}V(2,\,0,\, a,\, b,\, x)+4bV(1,\,1,\, a,\, b,\, x)+\frac{8}{b}V(0,2,\, a,\, b,\, x)\end{array}\label{eq:int_Gamma5/x^2_alt_anyl}\end{equation}

Note that the second and third terms in the expression on the right-hand
side differ formally from the first-through-third terms in the step-up
method (\ref{eq:int_Gamma5/x^2_is_anyl}). But expanding both versions
into their explicit forms and rearranging shows that they are indeed
equal:

\begin{equation}
\begin{array}{ccc}
4bV(1,\,1,\, a,\, b,\, x)+\frac{8}{b}V(0,2,\, a,\, b,\, x)\hspace{6.4cm}\\
=4b^{2}x\left(\frac{1}{bx}+1\right)e^{-ax^{2}-bx}+8b\left(\frac{3}{b^{2}x^{2}}+\frac{3}{bx}+1\right)e^{-ax^{2}-bx}\\ \\
12bV(0,\,0,\, a,\, b,\, x)+\frac{24}{b}V(-1,\,1,\, a,\, b,\, x)+4b^{2}V(1,\,0,\, a,\, b,\, x)\hspace{2.4cm}\\
\hspace{-0.7cm}=4b^{2}xe^{-ax^{2}-bx}+12be^{-ax^{2}-bx}+\frac{24\left(\frac{1}{bx}+1\right)e^{-ax^{2}-bx}}{x}\end{array}\label{eq:int_Gamma5/x^2_alt_anyl_comparsion}\end{equation}

One may continue on in this fashion via any of these three approaches
to develop a set of indefinite integrals for the incomplete gamma
function for indices and powers as large as one likes. We turn, instead,
to other Bessel functions.

\section{Modified Bessel Functions }

We move now from half-integer Macdonald functions (modified spherical
Bessel function of the second kind) to the indefinite integrals of
modified spherical Bessel function of the first kind (again given
in each of three equivalent forms): %
{}

\subsection*{Corollary 29}

\begin{equation}
\begin{array}{ll}
\sqrt{2\pi}\int x^{\frac{1}{2}}e^{-ax^{2}}b^{\frac{1}{2}}I_{\frac{1}{2}}(xb)\, dx\\
=\sqrt{2\pi}\int i^{-\frac{1}{2}}\sqrt{b}\sqrt{x}e^{-ax^{2}}G_{0,2}^{1,0}\left(-\frac{1}{4}b^{2}x^{2}|\begin{array}{c}
\frac{1}{4},-\frac{1}{4}\end{array}\right)\, dx\\
=\int2e^{-ax^{2}}\sinh(bx)\, dx=\frac{\sqrt{\pi}e^{\frac{b^{2}}{4a}}}{2\sqrt{a}}\left(\text{erf}\left(\frac{2ax-b}{2\sqrt{a}}\right)-\text{erf}\left(\frac{2ax+b}{2\sqrt{a}}\right)\right)\\
=V(0,\,0,\, a,\,-b,\, x)-V(0,\,0,\, a,\, b,\, x)\end{array}\:.\label{eq:I[1/2] x^1/2}\end{equation}

\noindent The final step in the fourth line follows from the definition~
\cite{GR5} (p. 29 No. 1.311.2)

\begin{equation}
\sinh(bx)=\frac{1}{2}\left(e^{bx}-e^{-bx}\right)\label{eq:sinh}\end{equation}

\noindent and the fact that integrals with any non-negative integer
power of \emph{x} multiplying $e^{-ax^{2}-bx}$ in Prudnikov, Brychkov,
and Marichev, \cite{PBM1} (p. 140 No. 1.3.3.16,17) have $a>0$ but
no such restriction on \emph{b}. For the same reason

\subsection*{Corollary 30}

\begin{equation}
\begin{array}{ll}
\sqrt{2\pi}\int x^{\frac{1}{2}+1}e^{-ax^{2}}b^{\frac{1}{2}}I_{\frac{1}{2}}(xb)\, dx\\
=\sqrt{2\pi}\int i^{-\frac{1}{2}}\sqrt{b}x^{\frac{1}{2}+1}e^{-ax^{2}}G_{0,2}^{1,0}\left(-\frac{1}{4}b^{2}x^{2}|\begin{array}{c}
\frac{1}{4},-\frac{1}{4}\end{array}\right)\, dx\\
=\int2e^{-ax^{2}}x\,\sinh(bx)\, dx=V(1,\,0,\, a,\,-b,\, x)-V(1,\,0,\, a,\, b,\, x)\end{array}\:\label{eq:I[1/2] x^1/2+1}\end{equation}

\noindent and

\subsection*{Corollary 31}

\begin{equation}
\begin{array}{ll}
\sqrt{2\pi}\int x^{\frac{1}{2}+2}e^{-ax^{2}}b^{\frac{1}{2}}I_{\frac{1}{2}}(xb)\, dx\\
=\sqrt{2\pi}\int i^{-\frac{1}{2}}\sqrt{b}x^{\frac{1}{2}+2}e^{-ax^{2}}G_{0,2}^{1,0}\left(-\frac{1}{4}b^{2}x^{2}|\begin{array}{c}
\frac{1}{4},-\frac{1}{4}\end{array}\right)\, dx\\
=\int2e^{-ax^{2}}x^{2}\,\sinh(bx)\, dx=V(2,\,0,\, a,\,-b,\, x)-V(2,\,0,\, a,\, b,\, x)\end{array}\:.\label{eq:I[1/2] x^1/2+2}\end{equation}

Similarly, one may find

\subsection*{Corollary 32}

\begin{equation}
\begin{array}{ll}
\sqrt{2\pi}\int x^{\frac{1}{2}-1}e^{-ax^{2}}b^{\frac{3}{2}}I_{\frac{3}{2}}(xb)\, dx\\
=\sqrt{2\pi}\int i^{-\frac{1}{2}}x^{\frac{1}{2}-1}e^{-ax^{2}}b^{\frac{3}{2}}G_{0,2}^{1,0}\left(-\frac{1}{4}b^{2}x^{2}|\begin{array}{c}
\frac{3}{4},-\frac{3}{4}\end{array}\right)\, dx\\
=\int\frac{1}{x}be^{-ax^{2}}\left(2\cosh(bx)-\frac{2\sinh(bx)}{bx}\right)\, dx=\frac{\sqrt{\pi}\sqrt{a}xe^{\frac{b^{2}}{4a}}\text{erf}\left(\frac{2ax-b}{2\sqrt{a}}\right)+e^{bx-ax^{2}}}{x}-\frac{\sqrt{\pi}\sqrt{a}xe^{\frac{b^{2}}{4a}}\text{erf}\left(\frac{2ax+b}{2\sqrt{a}}\right)+e^{-ax^{2}-bx}}{x}\\
=V(-1,\,1,\, a,\,-b,\, x)-V(-1,\,1,\, a,\, b,\, x)\end{array}\:,\label{eq:I[3/2] x^(1/2-1)}\end{equation}

\noindent and likewise for the full set.

On the other hand, the in-between integrals for $I_{\frac{n}{2}}(xb)$
are most easily determined by examining the form of the integrand
involving hyperbolic functions and using the even and odd series expansions
defined in the last two lines of (\ref{eq:int exp-ax^2 series of cosh-sinh}),
as in

\subsection*{Corollary 33}

\begin{equation}
\begin{array}{ll}
\sqrt{2\pi}\int x^{\frac{1}{2}}e^{-ax^{2}}b^{\frac{3}{2}}I_{\frac{3}{2}}(xb)\, dx\\
=\sqrt{2\pi}\int i^{-\frac{1}{2}}x^{\frac{1}{2}}e^{-ax^{2}}b^{\frac{3}{2}}G_{0,2}^{1,0}\left(-\frac{1}{4}b^{2}x^{2}|\begin{array}{c}
\frac{3}{4},-\frac{3}{4}\end{array}\right)\, dx\\
=\int be^{-ax^{2}}\left(2\cosh(bx)-\frac{2\sinh(bx)}{bx}\right)\, dx\\
=2b\Upsilon_{e}\left(0,\, a,\, b,\, x\right)+2\Upsilon_{o}\left(1,\, a,\, b,\, x\right)\end{array}\:.\label{eq:I[3/2] x^(1/2-0)}\end{equation}

\section{Bessel Functions }

Indefinite integrals of spherical Bessel functions of the first kind
follow from their modified kin.

\subsection*{Corollary 34}

\begin{equation}
\begin{array}{ll}
\sqrt{2\pi}\int x^{\frac{1}{2}}e^{-ax^{2}}b^{\frac{1}{2}}J_{\frac{1}{2}}(xb)\, dx\\
=\sqrt{2\pi}\int\sqrt{b}\sqrt{x}e^{-ax^{2}}G_{0,2}^{1,0}\left(\frac{1}{4}b^{2}x^{2}|\begin{array}{c}
\frac{1}{4},-\frac{1}{4}\end{array}\right)\, dx\\
=\int2e^{-ax^{2}}\sin(bx)\, dx=-\frac{\sqrt{\pi}}{2\sqrt{a}}e^{-\frac{b^{2}}{4a}}\left(\text{erfi}\left(\frac{b+2iax}{2\sqrt{a}}\right)-i\text{erf}\left(\frac{2ax+ib}{2\sqrt{a}}\right)\right)\\
=i\frac{\sqrt{\pi}}{2\sqrt{a}}e^{-\frac{b^{2}}{4a}}\left(\text{erf}\left(\frac{2ax+ib}{2\sqrt{a}}\right)-\text{erf}\left(\frac{2ax-ib}{2\sqrt{a}}\right)\right)\\
=iV(0,\,0,\, a,\, ib,\, x)-iV(0,\,0,\, a,\,-ib,\, x)\end{array}\:.\label{eq:J[1/2] x^1/2}\end{equation}

\noindent The fourth line follows from the conversions \cite{erfi}

\begin{equation}
\text{erfi}(z)=-i\,\text{erf}(iz)\label{eq:erfi}\end{equation}

\noindent and \cite{erf(-z)}

\begin{equation}
\text{erf}(-z)=-\text{erf}(z)\label{eq:erf(-z)}\end{equation}

\noindent and the last line can simply be read off the fourth in comparison
with the last two lines of (\ref{eq:I[1/2] x^1/2}). Equivalently,
one may obtain it from the third line and the definition \cite{GR5}
(p. 29 No. 1.311.1)

\begin{equation}
\sin(z)=-i\sinh(iz)\label{eq:sin_from_sinh}\end{equation}

\noindent in (\ref{eq:I[1/2] x^1/2}), while noting that integrals
with any non-negative integer power of \emph{x} multiplying $e^{-ax^{2}-bx}$
in Prudnikov, Brychkov, and Marichev, \cite{PBM1}(p. 140 No. 1.3.3.16,17)
require \emph{a} to be real (and positive) but there is no such restriction
on \emph{b}.

Although the form on the fourth line of (\ref{eq:J[1/2] x^1/2}),
and other integrals that follow, is convenient for tying new integrals
to prior ones (the final line), a more compact form may be found by
using \cite{Im[Erf(x+iy)]}

\begin{equation}
\Im(\text{erf}(x+iy))=\frac{1}{2}i(\text{erf}(x-iy)-\text{erf}(x+iy))\label{eq:Im[Erf(x+iy)]}\end{equation}

\noindent so that

\begin{equation}
i\left(\text{erf}\left(\frac{2ax+ib}{2\sqrt{a}}\right)-\text{erf}\left(\frac{2ax-ib}{2\sqrt{a}}\right)\right)=-2\Im\left(\text{erf}\left(\frac{2ax+ib}{2\sqrt{a}}\right)\right)\:.\label{eq:Im[erf[a+ib]]}\end{equation}

\noindent After letting $b\rightarrow-b$, this casts (\ref{eq:J[1/2] x^1/2})
into a form that is a generalization of the upper line of the only
tabled integral of the class presented in the present work \cite{PBM1}
(p. 234 No. 1.5.49.17):

\begin{equation}
\begin{array}{ll}
\int e^{-x^{2}}\left\{ \begin{array}{c}
\sin(bx)\\
\cos(bx)\end{array}\right\} \, dx=\frac{\sqrt{\pi}}{2}e^{-\frac{b^{2}}{4}}\left\{ \begin{array}{c}
\Im\left(\text{erf}\left(x-\frac{ib}{2}\right)\right)\\
\Re\left(\text{erf}\left(x-\frac{ib}{2}\right)\right)\end{array}\right\} \end{array}\:.\label{eq:int exp[-ax^2](sin or cos)}\end{equation}

Following the pattern in (\ref{eq:J[1/2] x^1/2})

\subsection*{Corollary 35}

\begin{equation}
\begin{array}{ll}
\sqrt{2\pi}\int x^{\frac{1}{2}+1}e^{-ax^{2}}b^{\frac{1}{2}}J_{\frac{1}{2}}(xb)\, dx\\
=\sqrt{2\pi}\int\sqrt{b}x^{\frac{1}{2}+1}e^{-ax^{2}}G_{0,2}^{1,0}\left(\frac{1}{4}b^{2}x^{2}|\begin{array}{c}
\frac{1}{4},-\frac{1}{4}\end{array}\right)\,\, dx\\
=\int2e^{-ax^{2}}x\,\sin(bx)\, dx=iV(1,\,0,\, a,\, ib,\, x)-iV(1,\,0,\, a,\,-ib,\, x)\end{array}\:\label{eq:J[1/2] x^1/2+1}\end{equation}

\noindent and

\subsection*{Corollary 36}

\begin{equation}
\begin{array}{ll}
\sqrt{2\pi}\int x^{\frac{1}{2}+2}e^{-ax^{2}}b^{\frac{1}{2}}J_{\frac{1}{2}}(xb)\, dx\\
=\sqrt{2\pi}\int\sqrt{b}x^{\frac{1}{2}+2}e^{-ax^{2}}G_{0,2}^{1,0}\left(\frac{1}{4}b^{2}x^{2}|\begin{array}{c}
\frac{1}{4},-\frac{1}{4}\end{array}\right)\, dx\\
=\int2e^{-ax^{2}}x^{2}\,\sin(bx)\, dx=iV(2,\,0,\, a,\, ib,\, x)-iV(2,\,0,\, a,\,-ib,\, x)\end{array}\:.\label{eq:J[1/2] x^1/2+2}\end{equation}

Similarly, one finds

\subsection*{Corollary 37}

\begin{equation}
\begin{array}{ll}
\sqrt{2\pi}\int x^{\frac{1}{2}-1}e^{-ax^{2}}b^{\frac{3}{2}}J_{\frac{3}{2}}(xb)\, dx\\
=\sqrt{2\pi}\int x^{\frac{1}{2}-1}e^{-ax^{2}}b^{\frac{3}{2}}G_{0,2}^{1,0}\left(\frac{1}{4}b^{2}x^{2}|\begin{array}{c}
\frac{3}{4},-\frac{3}{4}\end{array}\right)\, dx\\
=\int\frac{1}{x}be^{-ax^{2}}\left(-2\cos(bx)+\frac{2\sin(bx)}{bx}\right)\, dx\\
=i\sqrt{\pi}\sqrt{a}e^{-\frac{b^{2}}{4a}}\text{erf}\left(\frac{2ax-ib}{2\sqrt{a}}\right)-i\sqrt{\pi}\sqrt{a}e^{-\frac{b^{2}}{4a}}\text{erf}\left(\frac{2ax+ib}{2\sqrt{a}}\right)+\frac{i}{x}\left(e^{-ax^{2}+ibx}-e^{-ax^{2}-ibx}\right)\\
=iV(-1,\,1,\, a,\, ib,\, x)-iV(-1,\,1,\, a,\,-ib,\, x)\end{array}\:,\label{eq:J[3/2] x^(1/2-1)}\end{equation}

\noindent and likewise for the full set.

Again, the in-between integrals for $J_{\frac{n}{2}}(xb)$ are most
easily determined by examining the form of the integrand involving
trigonometric functions, using the relations \cite{GR5} (p. 29 No.
1.311.1,3\begin{equation}
\begin{array}{ccc}
\sin(bx) & = & -i\sinh(ibx)\\
\cos(bx) & = & \cosh(ibx)\end{array}\label{eq:cos2cosh}\end{equation}

\noindent and using the even and odd series expansions defined in
the last two lines of (\ref{eq:int exp-ax^2 series of cosh-sinh}), as in

\subsection*{Corollary 38}

\begin{equation}
\begin{array}{ll}
\sqrt{2\pi}\int x^{\frac{1}{2}}e^{-ax^{2}}b^{\frac{3}{2}}J_{\frac{3}{2}}(xb)\, dx\\
=\sqrt{2\pi}\int x^{\frac{1}{2}}e^{-ax^{2}}b^{\frac{3}{2}}G_{0,2}^{1,0}\left(\frac{1}{4}b^{2}x^{2}|\begin{array}{c}
\frac{3}{4},-\frac{3}{4}\end{array}\right)\, dx\\
=\int be^{-ax^{2}}\left(-2\cos(bx)+\frac{2\sin(bx)}{bx}\right)\, dx\\
=\int be^{-ax^{2}}\left(-2\cosh(ibx)-\frac{2i\sinh(ibx)}{bx}\right)\, dx\\
=-2b\Upsilon_{e}\left(0,\, a,\, ib,\, x\right)+2i\Upsilon_{o}\left(1,\, a,\, ib,\, x\right)\end{array}\:.\label{eq:J[3/2] x^(1/2-0)}\end{equation}

\section{Neumann Functions }

Indefinite integrals of spherical Bessel functions of the second kind
$Y_{\frac{n}{2}}(xb)=N_{\frac{n}{2}}(xb)$ parallel the prior section.

\subsection*{Corollary 39}

\begin{equation}
\begin{array}{ll}
\sqrt{2\pi}\int x^{\frac{1}{2}}e^{-ax^{2}}b^{\frac{1}{2}}Y_{\frac{1}{2}}(xb)\, dx\\
=\sqrt{2\pi}\int\sqrt{b}\sqrt{x}e^{-ax^{2}}G_{1,3}^{2,0}\left(\frac{b^{2}x^{2}}{4}|\begin{array}{c}
-\frac{3}{4}\\
\frac{1}{4},-\frac{1}{4},-\frac{3}{4}\end{array}\right)\, dx\\
=-\int2e^{-ax^{2}}\cos(bx)\, dx=-\frac{\sqrt{\pi}e^{-\frac{b^{2}}{4a}}\left(\text{erf}\left(\frac{2ax-ib}{2\sqrt{a}}\right)+i\text{erfi}\left(\frac{b-2iax}{2\sqrt{a}}\right)\right)}{2\sqrt{a}}\\
=-\frac{\sqrt{\pi}}{2\sqrt{a}}e^{-\frac{b^{2}}{4a}}\left(\text{erf}\left(\frac{2ax-ib}{2\sqrt{a}}\right)+\text{erf}\left(\frac{2ax+ib}{2\sqrt{a}}\right)\right)\\
=-V(0,\,0,\, a,\, ib,\, x)-V(0,\,0,\, a,\,-ib,\, x)\end{array}\:.\label{eq:Y[1/2] x^1/2}\end{equation}

\noindent A compact form may be found by using \cite{Re[Erf(x+iy)]}\begin{equation}
\Re(\text{erf}(x+iy))=\frac{1}{2}(\text{erf}(x+iy)+\text{erf}(x-iy))\label{eq:Re[Erf(x+iy)]}\end{equation}
 so that\begin{equation}
\text{erf}\left(\frac{2ax-ib}{2\sqrt{a}}\right)+\text{erf}\left(\frac{2ax+ib}{2\sqrt{a}}\right)=2\Re\left(\text{erf}\left(\frac{2ax+ib}{2\sqrt{a}}\right)\right)\:.\label{eq:Im[erf[a+ib]]-1}\end{equation}

\noindent Since $\cos(-bx)=\cos(bx)$, one may let $b\rightarrow-b$
throughout to obtain the lower line of (\ref{eq:int exp[-ax^2](sin or cos)})
when $a=1$.

The process leading to (\ref{eq:Y[1/2] x^1/2}) likewise gives

\subsection*{Corollary 40}

\begin{equation}
\begin{array}{ll}
\sqrt{2\pi}\int x^{\frac{1}{2}+1}e^{-ax^{2}}b^{\frac{1}{2}}Y_{\frac{1}{2}}(xb)\, dx\\
=\sqrt{2\pi}\int\sqrt{b}x^{\frac{1}{2}+1}e^{-ax^{2}}G_{1,3}^{2,0}\left(\frac{b^{2}x^{2}}{4}|\begin{array}{c}
-\frac{3}{4}\\
\frac{1}{4},-\frac{1}{4},-\frac{3}{4}\end{array}\right)\,\, dx\\
=-\int2e^{-ax^{2}}x\,\cos(bx)\, dx=-V(1,\,0,\, a,\, ib,\, x)-V(1,\,0,\, a,\,-ib,\, x)\end{array}\:,\label{eq:Y[1/2] x^1/2+1}\end{equation}

\subsection*{Corollary 41}

\begin{equation}
\begin{array}{ll}
\sqrt{2\pi}\int x^{\frac{1}{2}+2}e^{-ax^{2}}b^{\frac{1}{2}}Y_{\frac{1}{2}}(xb)\, dx\\
=\sqrt{2\pi}\int\sqrt{b}x^{\frac{1}{2}+2}e^{-ax^{2}}G_{1,3}^{2,0}\left(\frac{b^{2}x^{2}}{4}|\begin{array}{c}
-\frac{3}{4}\\
\frac{1}{4},-\frac{1}{4},-\frac{3}{4}\end{array}\right)\, dx\\
=-\int2e^{-ax^{2}}x^{2}\,\cos(bx)\, dx=-V(2,\,0,\, a,\, ib,\, x)-V(2,\,0,\, a,\,-ib,\, x)\end{array}\:,\label{eq:Y[1/2] x^1/2+2}\end{equation}

\subsection*{Corollary 42}

\begin{equation}
\begin{array}{ll}
\sqrt{2\pi}\int x^{\frac{1}{2}-1}e^{-ax^{2}}b^{\frac{3}{2}}Y_{\frac{3}{2}}(xb)\, dx\\
=\sqrt{2\pi}\int x^{\frac{1}{2}-1}e^{-ax^{2}}b^{\frac{3}{2}}G_{1,3}^{2,0}\left(\frac{b^{2}x^{2}}{4}|\begin{array}{c}
-\frac{5}{4}\\
\frac{3}{4},-\frac{3}{4},-\frac{5}{4}\end{array}\right)\, dx\\
=\int\frac{1}{x}be^{-ax^{2}}\left(-2\sin(bx)-2\frac{\cos(bx)}{bx}\right)\, dx\\
=\sqrt{\pi}\sqrt{a}e^{-\frac{b^{2}}{4a}}\text{eErf}\left(\frac{2ax-ib}{2\sqrt{a}}\right)+\sqrt{\pi}\sqrt{a}e^{-\frac{b^{2}}{4a}}\text{erf}\left(\frac{2ax+ib}{2\sqrt{a}}\right)+\frac{1}{x}\left(e^{-ax^{2}+(-i)bx}+e^{-ax^{2}+ibx}\right)\\
=-V(-1,\,1,\, a,\, ib,\, x)-V(-1,\,1,\, a,\,-ib,\, x)\end{array}\:,\label{eq:Y[3/2] x^(1/2-1)}\end{equation}

\noindent and so on for the full set.

As with the Bessel functions of the first kind, the in-between integrals
for $Y_{\frac{n}{2}}(xb)$ are most easily determined by converting
from trigonometric to hyperbolic functions and examining the form
of the integrand in comparison to the even and odd series expansions
defined in the last two lines of (\ref{eq:int exp-ax^2 series of cosh-sinh}
):

\subsection*{Corollary 43}

\begin{equation}
\begin{array}{ll}
\sqrt{2\pi}\int x^{\frac{1}{2}}e^{-ax^{2}}b^{\frac{3}{2}}Y_{\frac{3}{2}}(xb)\, dx\\
=\sqrt{2\pi}\int x^{\frac{1}{2}}e^{-ax^{2}}b^{\frac{3}{2}}G_{1,3}^{2,0}\left(\frac{b^{2}x^{2}}{4}|\begin{array}{c}
-\frac{5}{4}\\
\frac{3}{4},-\frac{3}{4},-\frac{5}{4}\end{array}\right)\, dx\\
=\int be^{-ax^{2}}\left(-\frac{2\cos(bx)}{bx}-2\sin(bx)\right)\, dx\\
=\int be^{-ax^{2}}\left(-\frac{2\cosh(ibx)}{bx}+2i\sinh(ibx)\right)\, dx\\
=-2\Upsilon_{e}\left(1,\, a,\, ib,\, x\right)-2ib\Upsilon_{o}\left(0,\, a,\, ib,\, x\right)\end{array}\:.\label{eq:Y[3/2] x^(1/2-0)-1}\end{equation}

\section{Struve Functions }

Indefinite integrals of half-integer Struve H functions can be found
from those of Neumann Functions since \cite{GR5} (p. 997 No .8.552.3)\begin{equation}
\mathbf{H}_{n+\frac{1}{2}}(bx)=Y_{n+\frac{1}{2}}(bx)+\frac{1}{\pi}\sum_{m=0}^{n}\frac{2^{2m-n+\frac{1}{2}}\Gamma\left(m+\frac{1}{2}\right)b^{-2m+n-\frac{1}{2}}}{\Gamma(-m+n+1)}x^{-2m+n-\frac{1}{2}}\label{eq:H_as_Y+sum}\end{equation}

\noindent so that

\subsection*{Corollary 44}

\noindent \begin{equation}
\begin{array}{ll}
\sqrt{2\pi}\int x^{\frac{1}{2}}e^{-ax^{2}}b^{\frac{1}{2}}\mathbf{H}_{\frac{1}{2}}(bx)\, dx\\
=\sqrt{2\pi}\int\sqrt{b}\sqrt{x}e^{-ax^{2}}G_{1,3}^{1,1}\left(\frac{b^{2}x^{2}}{4}|\begin{array}{c}
\frac{3}{4}\\
\frac{3}{4},-\frac{1}{4},\frac{1}{4}\end{array}\right)\, dx\\
=-\int2e^{-ax^{2}}(\cos(bx)-1)\, dx\\
=-\frac{\sqrt{\pi}e^{-\frac{b^{2}}{4a}}\text{erf}\left(\frac{2ax+ib}{2\sqrt{a}}\right)}{2\sqrt{a}}+\frac{i\sqrt{\pi}e^{-\frac{b^{2}}{4a}}\text{erfi}\left(\frac{b+2iax}{2\sqrt{a}}\right)}{2\sqrt{a}}+\frac{\sqrt{\pi}\text{erf}\left(\sqrt{a}x\right)}{\sqrt{a}}\\
=-\frac{\sqrt{\pi}e^{-\frac{b^{2}}{4a}}\text{erf}\left(\frac{2ax+ib}{2\sqrt{a}}\right)}{2\sqrt{a}}-\frac{\sqrt{\pi}e^{-\frac{b^{2}}{4a}}\text{erf}\left(\frac{2ax-ib}{2\sqrt{a}}\right)}{2\sqrt{a}}+\frac{\sqrt{\pi}\text{erf}\left(\sqrt{a}x\right)}{\sqrt{a}}\\
=-V(0,\,0,\, a,\, ib,\, x)-V(0,\,0,\, a,\,-ib,\, x)+2V(0,\,0,\, a,\,0,\, x)\end{array}\:.\label{eq:H[1/2] x^1/2}\end{equation}

For $n>0$ one needs \cite{PBM1} (p. 140 No. 1.3.3.12 which I have
transformed from $e^{-a^{2}x^{2}}\rightarrow e^{-ax^{2}}$)\begin{equation}
\int\frac{e^{-ax^{2}}}{x^{2q+1}}\, dx=\frac{1}{2aq!}e^{-ax^{2}}\sum_{k=1}^{q}(-1)^{k}a^{k}(q-k)!x^{2k-2q-2}+\frac{\left((-1)^{q}a^{q}\right)\text{Ei}\left(-ax^{2}\right)}{2q!}\label{eq:exp[-a x^2]/x^{2 q+1}}\end{equation}

\noindent in which one simply omits the first term if $q=0$.

Indefinite integrals of half-integer Struve L functions can be found
from half-integer modified Bessel functions since~\cite{MagnusOberhettingerSoni}\begin{equation}
\mathbf{L}_{n+\frac{1}{2}}(bx)=I_{n+\frac{1}{2}}(bx)+\frac{2(-1)^{n}}{\pi}K_{n+\frac{1}{2}}(bx)-\frac{1}{\sqrt{\pi}}\sum_{m=0}^{n}\frac{(-1)^{m}\left(2^{-2m}(2m)!\right)2^{2m-n+\frac{1}{2}}b^{-2m+n-\frac{1}{2}}x^{-2m+n-\frac{1}{2}}}{m!(n-m)!}\label{eq:L as I K series}\end{equation}

\noindent so that

\subsection*{Corollary 45}

\noindent \begin{equation}
\begin{array}{ll}
\sqrt{2\pi}\int x^{\frac{1}{2}}e^{-ax^{2}}b^{\frac{1}{2}}\mathbf{L}_{\frac{1}{2}}(bx)\, dx\\
=-\int2\pi^{3/2}\sqrt{b}\sqrt{x}e^{-ax^{2}}G_{2,4}^{1,1}\left(\frac{b^{2}x^{2}}{4}|\begin{array}{c}
\frac{3}{4},\frac{1}{2}\\
\frac{3}{4},\frac{1}{2},-\frac{1}{4},\frac{1}{4}\end{array}\right)\, dx\\
=\int\left(2e^{-ax^{2}}\sinh(bx)+2e^{-ax^{2}-bx}-2e^{-ax^{2}}\right)\, dx\\
=\int\left(2e^{-ax^{2}}\cosh(bx)-2e^{-ax^{2}}\right)\, dx\\
=\frac{\sqrt{\pi}e^{\frac{b^{2}}{4a}}\text{erf}\left(\frac{2ax+b}{2\sqrt{a}}\right)}{2\sqrt{a}}+\frac{\sqrt{\pi}e^{\frac{b^{2}}{4a}}\text{erf}\left(\frac{2ax-b}{2\sqrt{a}}\right)}{2\sqrt{a}}-\frac{\sqrt{\pi}\text{erf}\left(\sqrt{a}x\right)}{\sqrt{a}}\\
=V(0,\,0,\, a,\, b,\, x)+V(0,\,0,\, a,\,-b,\, x)-2V(0,\,0,\, a,\,0,\, x)\end{array}\:.\label{eq:L[1/2] x^1/2}\end{equation}

\section*{Conclusion}

We have crafted a set of indefinite integrals for 63 half-integer
Bessel 
and Struve 
functions, and incomplete gamma 
functions with integer indices, each multiplied by $Exp(-ax^{2})$
and divided by powers. A series solution is given for the individual
terms (of any power) in such functions, split into even and odd portions,
which converges faster than conventional series derived from expanding
$Exp(-ax^{2})$ or $Exp(-bx)$. Eight integrals involving these series
are given. Three infinite sets of non-Bessel-related auxiliary indefinite
integrals are given in Appendix B, along with two cases given explicitly,
as are two indefinite integrals not related to those infinite sets.

The integrals of the sort presented in this paper do not lend themselves
to solution using computer calculus systems. Indeed, I could not get
Maple 2026 or Maxima to do even the most elementary of these integrals.
I was, however, able to carefully construct symbolic integration methods
using Mathematica, and have laid out the necessary input forms in
Section 6 in detail. This is both because one always hopes to have
colleagues take one's work onto unexpected paths. but also because
one needs cross-checking for one's own work.

\vspace{6pt}

\noindent \textbf{Funding:} This research received no external funding.

\noindent \textbf{Conflicts of interest: }The author declares no conflicts
of interest.

\section*{Appendix A. Mathematica code to calculate $\int x^{n}e^{-ax^{2}-bx}\, dx$}

Here we provide the Mathematica code to calculate $\int x^{n}e^{-ax^{2}-bx-c}\, dx$
for any non-negative integer \emph{n} from Korotkov's expression \cite{Korotkov}
(p. 4 No. 1.2.4 with $a\rightarrow\alpha^{2}$ and $b\rightarrow-\beta$)).
One is cautioned, however, that Mathematica 7 and 13 both have the
peculiarity of doubling the 

\begin{verbatim}

Korotkov124a[ n_, a_, b_,c_, x_] :=
-2^-n E^(-c - b x - a x^2) (
   Table[((n! (K - 1)! (-b)^(n + 1 - 2 K)) 
      Sum[ 
                 (4^(K - L) (2 x a + b)^(2 L - 2))/((L - 1)! (a)^(n - K + L))
        , {L, 1, K}])/((2 K - 1)! (n + 1 - 2 K)!)
    , {K, 1, n - Floor[n/2]}] + 
   Table[((n! (-b)^(n - 2 K))  
      Sum[
               (L! (2 x a + b)^(2 L - 1))/((2 L)! (a)^(n - K + L))
       , {L, 1,  K}])/( K! (n - 2 K)!)
    , {K, 1, Floor[n/2]}])

\end{verbatim}

\noindent
portion of the full expression whenever $n\geq3$ if it is evaluated
at the same time as the

\begin{verbatim}

Korotkov124b[n_, a_, b_, c_, x_] :=
1/Sqrt[a] 2^(-1 - n) E^(b^2/(4 a) - c) Sqrt[\[Pi]] Erf[(b + 2 a x)/(2 Sqrt[a])] *
Table[(n! (-b)^(n - 2 K))/(
  K! (n - 2 K)! a^(n - K))
, {K, 0, 1}]

\end{verbatim}

\noindent
portion. If these portions are evaluated separately for a given \emph{n}
and then added, one gets the correct result. I have retain the constant
prefeactor $e^{-c}$ from Korotkov's expression but set $c=0$ in
my own work.

The function {}``Table'' is used at several places in the code because
if one uses the function {}``Sum'' instead, Mathematica will try
to calculate the analytical result of summing the series. For the
simplest case of $n=1$, using Table takes 3.2 milliseconds whereas
using Sum takes six minutes. and also givs a set of incomplete gamma
functions instead of the desired error functions. For $n=2$, for
instance, {}``Sum'' gives 

\begin{equation}
\begin{array}{l}
-\frac{1}{4}e^{-ax^{2}-bx-c}\left(\frac{e^{\frac{(2ax+b)^{2}}{4a}}\left(24\sqrt{\pi}a^{3/2}\left(\frac{(2ax+b)^{2}}{a}\right)^{3/2}\text{erf}\left(\frac{2ax+b}{2\sqrt{a}}\right)-48(2ax+b)^{3}\left(\frac{\sqrt{\pi}}{2}-\Gamma\left(\frac{3}{2},\frac{(b+2ax)^{2}}{4a}\right)\right)\right)}{24a^{3}\left(\frac{(2ax+b)^{2}}{a}\right)^{3/2}}-\frac{2b}{a^{2}}\right)\end{array}\:,\label{eq:int x^2 exp as gamma 1st part}\end{equation}
for the first part of the code whose second term contains the function
$\Gamma\left(\frac{3}{2},\frac{(b+2ax)^{2}}{4a}\right)$ that must
be converted to $\text{erf}\left(\frac{2ax+b}{2\sqrt{a}}\right)$
by hand \cite{gamma to erfc,erfc to erf} to cancel the first term.

On the other hand, {}``Table'' gives

\begin{equation}
\left\{ -\frac{1}{4}\left(\frac{2ax+b}{a^{2}}-\frac{2b}{a^{2}}\right)e^{-ax^{2}-bx-c}\right\} \:.\label{eq:int x^2 exp using Table}\end{equation}
This single-entry list (indicated by the curly brackets) is turned
into a sum with the {}``Total{[}\%{]}'' command, where {}``\%''
refers to the previous result. This is added to the second part of
the code

\begin{equation}
\frac{\sqrt{\pi}b^{2}\left(\frac{a}{b^{2}}+\frac{1}{2}\right)e^{\frac{b^{2}}{4a}-c}\text{erf}\left(\frac{2ax+b}{2\sqrt{a}}\right)}{4a^{5/2}}\label{eq:2nd part n=00003D2}\end{equation}
to give the final result.

For $n>2$ there will be internal tables (indicated by curly brackets)
that must be totalled, before totaling the external one. E.g, for
$n=3$ the first portion of the code gives

\begin{quote}
-(1/8) E\textasciicircum{}(-c - b x - a x\textasciicircum{}2) (\{-((3
b (b + 2 a x))/a\textasciicircum{}3)\} + \{(3 b\textasciicircum{}2)/a\textasciicircum{}3,
4/a\textasciicircum{}2 + (b + 2 a x)\textasciicircum{}2/a\textasciicircum{}3\})
\end{quote}

\noindent
within which one inserts the Total command as

\begin{quote}
-(1/8) E\textasciicircum{}(-c - b x - a x\textasciicircum{}2) (Total{[}\{-((3
b (b + 2 a x))/a\textasciicircum{}3)\}{]} + Total{[}\{(3 b\textasciicircum{}2)/a\textasciicircum{}3,
4/a\textasciicircum{}2 + (b + 2 a x)\textasciicircum{}2/a\textasciicircum{}3\}{]})
\end{quote}

\noindent
For $n>3$, three subsequent commands are useful for a polished result. 

\begin{quote}
Expand{[}\%{]}

\% /. \{E\textasciicircum{}(-x (b + a x)) -> E\textasciicircum{}(-b
x - a x\textasciicircum{}2), E\textasciicircum{}(b\textasciicircum{}2/(4
a) - (b + 2 a x)\textasciicircum{}2/(4 a)) -> E\textasciicircum{}(-b
x - a x\textasciicircum{}2), E\textasciicircum{}(-x (b + a x) + (b
+ 2 a x)\textasciicircum{}2/(4 a)) -> E\textasciicircum{}(b\textasciicircum{}2/(4
a)), E\textasciicircum{}(-b x - a x\textasciicircum{}2 + (b + 2 a
x)\textasciicircum{}2/(4 a)) -> E\textasciicircum{}(b\textasciicircum{}2/(4
a))\}

Collect{[}\%, E\textasciicircum{}(-b x - a x\textasciicircum{}2){]}
\end{quote}

\noindent
Use of the Table command, thus, creates the need for some additional
steps by hand, but this downside is far outweight by the faster result
in the desired form.

\section{Appendix B. Additional integrations by parts}

This appendix contains indefinite integrals not of the Bessel functions
themselves, but of related functions. The first set are those integrals
required to prove Theoems 6 and 10.

\subsection{Auxiliary indefinite integrals}

To prove Theorem 6 by first writing the error function in an integral
representation, then finding the negative of the indefinite integral
of the represented integrand with respect to \emph{a}, and then undoing
the integral representation by integrating that result respect to
\emph{x}, one needs this final integral. It may be shown to be: 

\textbf{Theorem 18}

\begin{equation}
\begin{array}{ll}
\int e^{-ax^{2}-bx}\left(-\frac{6}{x^{4}}-\frac{6a-b^{2}}{x^{2}}\right)\, dx\\
=e^{-ax^{2}-bx}\left(\frac{2a}{x}-\frac{b}{x^{2}}+\frac{2}{x^{3}}\right)+2\sqrt{\pi}a^{3/2}e^{\frac{b^{2}}{4a}}\text{erf}\left(\frac{2ax+b}{2\sqrt{a}}\right)\end{array}\:\label{eq:int for thm6 Theorem 18}\end{equation}

\textbf{Proof:}

One may apply the equations (\ref{eq:u}) through (\ref{eq:uv-intvdu}),
for integration by parts of Section 3, but with 

\begin{equation}
\text{dv}_{1}=\frac{6}{x^{4}}\:,\label{eq:dv for thm6 proof}\end{equation}
which reduces the inverse powers in the integral to

\begin{equation}
\frac{2e^{-ax^{2}-\text{bx}}}{x^{3}}+\int e^{-ax^{2}-bx}\left(\frac{4a}{x^{2}}+\frac{2b}{x^{3}}\right)\, dx-\int\frac{\left(6a-b^{2}\right)e^{-ax^{2}-bx}}{x^{2}}\, dx.\label{eq:int2 for thm6}\end{equation}
Then looping this process with 

\begin{equation}
\text{dv}_{2}=\frac{2b}{x^{3}}\:,\label{eq:dv2 for thm6 proof-1}\end{equation}
the further-reduced-power integral becomes 

\begin{equation}
\frac{2e^{-ax^{2}-\text{bx}}}{x^{3}}-\frac{be^{-ax^{2}-\text{bx}}}{x^{2}}-\int\frac{2ae^{-ax^{2}-bx}}{x^{2}}\, dx-\int\frac{2abe^{-ax^{2}-bx}}{x}\, dx\:.\label{eq:int3 for thm6-1}\end{equation}

Looping a second time with 

\begin{equation}
\text{dv}_{3}=\frac{2a}{x^{2}}\:,\label{eq:dv2 for thm6 proof}\end{equation}
gives, 

\begin{equation}
\frac{2e^{-ax^{2}-\text{bx}}}{x^{3}}-\frac{be^{-ax^{2}-\text{bx}}}{x^{2}}+\frac{2ae^{-ax^{2}-\text{bx}}}{x}+\int4a^{2}e^{-ax^{2}-bx}\, dx\:.\label{eq:int3 for thm6}\end{equation}
whose final integral is given by \cite{PBM1} (p. 140 No. 1.3.3.17)
(\ref{eq:PBM1 p. 140 No. 1.3.3.17}), completing the proof. $\square$

For Theorem 10, one needs to show the following: 

\textbf{Theorem 19}

\begin{equation}
\begin{array}{ll}
\int e^{-ax^{2}-bx}\left(\frac{40}{x^{6}}+\frac{40a-2b^{2}}{x^{4}}+\frac{20a^{2}-2ab^{2}}{x^{2}}\right)\, dx\\
=e^{-ax^{2}-bx}\left(-\frac{4a^{2}}{x}+\frac{2ab}{x^{2}}-\frac{8a}{x^{3}}+\frac{2b}{x^{4}}-\frac{8}{x^{5}}\right)-4\sqrt{\pi}a^{5/2}e^{\frac{b^{2}}{4a}}\text{erf}\left(\frac{2ax+b}{2\sqrt{a}}\right)\end{array}\:\label{eq:int for th10 Theorem 19}\end{equation}

\textbf{Proof:}

One may apply the equations (\ref{eq:u}) through (\ref{eq:uv-intvdu})
recursively. 

\begin{equation}
\begin{array}{ccc}
 &  & \int e^{-ax^{2}-bx}\left(\frac{40}{x^{6}}+\frac{40a-2b^{2}}{x^{4}}+\frac{20a^{2}-2ab^{2}}{x^{2}}\right)\, dx\\
 & = & -\frac{8e^{-ax^{2}-bx}}{x^{5}}-\int e^{-ax^{2}-bx}\left(\frac{8b}{x^{5}}+\frac{16a}{x^{4}}\right)\, dx+\int\frac{\left(20a^{2}-2ab^{2}\right)e^{-ax^{2}-bx}}{x^{2}}\, dx+\int\frac{\left(40a-2b^{2}\right)e^{-ax^{2}-bx}}{x^{4}}\, dx\\
 & = & -\frac{8e^{-ax^{2}-bx}}{x^{5}}+\frac{2be^{-ax^{2}-bx}}{x^{4}}+\int\frac{24ae^{-x(ax+b)}}{x^{4}}\, dx+\int\frac{4abe^{-ax^{2}-bx}}{x^{3}}\, dx+\int\frac{\left(20a^{2}-2ab^{2}\right)e^{-ax^{2}-bx}}{x^{2}}\, dx\\
 & = & -\frac{8e^{-ax^{2}-bx}}{x^{5}}+\frac{2be^{-ax^{2}-bx}}{x^{4}}-\frac{8ae^{-ax^{2}-bx}}{x^{3}}+\int-\frac{4abe^{-x(ax+b)}}{x^{3}}\, dx+\int\frac{2a\left(2a-b^{2}\right)e^{-x(ax+b)}}{x^{2}}\, dx\\
 & = & -\frac{8e^{-ax^{2}-bx}}{x^{5}}+\frac{2be^{-ax^{2}-bx}}{x^{4}}-\frac{8ae^{-ax^{2}-bx}}{x^{3}}+\frac{2abe^{-ax^{2}-bx}}{x^{2}}+\int\frac{4a^{2}e^{-x(ax+b)}}{x^{2}}\, dx+\int\frac{4a^{2}be^{-ax^{2}-bx}}{x}\, dx\\
 &  & -\frac{8e^{-ax^{2}-bx}}{x^{5}}+\frac{2be^{-ax^{2}-bx}}{x^{4}}-\frac{8ae^{-ax^{2}-bx}}{x^{3}}+\frac{2abe^{-ax^{2}-bx}}{x^{2}}-\frac{4a^{2}e^{-ax^{2}-bx}}{x}-\int8a^{3}e^{-ax^{2}-bx}\, dx\end{array}\label{eq:proof of thm 101}\end{equation}
whose final integral is given by \cite{PBM1} (p. 140 No. 1.3.3.17)
(\ref{eq:PBM1 p. 140 No. 1.3.3.17}), completing the proof. $\square$

\subsection{Other integrals of interest}

Integration by parts failed for Theorem 6, but one may use the seed
of that failure to establish several infinite sets of other integrals
that are not part of the context for this paper, but may be of use
to some readers who focus on different problems.

\textbf{Theorem 20}

\begin{equation}
\begin{array}{ll}
\int e^{-ax^{2}-bx}\left(x^{-\upsilon}-\frac{bx^{1-\upsilon}}{1-\upsilon}-\frac{2ax^{2-\upsilon}}{1-\upsilon}\right)\, dx\\
=\frac{x^{1-\upsilon}}{1-\upsilon}e^{-ax^{2}-bx}\quad\left[\upsilon\neq1\right]\end{array}\:\label{eq:generic int Theorem 20}\end{equation}

\textbf{Proof }

Consider again the integration by parts definitions of Section 2,

\begin{equation}
u=e^{-ax^{2}-bx}\label{eq:u-1}\end{equation}
so that \begin{equation}
\text{du}=e^{-ax^{2}-bx}(-2ax-b)\:.\label{eq:du-1}\end{equation}
Consider, then, the term having the most negative power within the
integrand, 

\begin{equation}
\text{dv}=x^{-\upsilon}\:,\label{eq:dv for generic}\end{equation}
whose integral is

\begin{equation}
v=\frac{x^{1-\upsilon}}{1-\upsilon}\quad\left[\upsilon\neq1\right].\label{eq:generic v}\end{equation}
Then the integral can be transformed to

\begin{equation}
\int u\, dv=uv-\int v\, du=\frac{x^{1-\upsilon}e^{-ax^{2}-bx}}{1-\upsilon}-\int e^{-ax^{2}-bx}\left(-\frac{bx^{1-\upsilon}}{1-\upsilon}-\frac{2ax^{2-\upsilon}}{1-\upsilon}\right)\, dx\:.\label{eq:uv-intvdu generic}\end{equation}
One may cancel the first term in the above integral by having a second
term in the original integral whose form is

\begin{equation}
-\int\frac{bx^{1-\upsilon}}{1-\upsilon}e^{-ax^{2}-bx}\, dx\:,\label{eq:second generic term}\end{equation}
and one may cancel the second term in (\ref{eq:uv-intvdu generic})
by having a third term in the original integral whose form is

\[
-\int\frac{2ax^{2-\upsilon}}{1-\upsilon}e^{-ax^{2}-bx}\, dx\:.\]
This completes the proof. $\square$

Several infininite series of integrals follow from Theorem 20 by integrating
with respect to paramaters. 

\textbf{Theorem 21}

\begin{equation}
\begin{array}{ll}
\int e^{-ax^{2}-bx}\left(\frac{\upsilon+1}{\upsilon-1}x^{-\upsilon-2}+\frac{bx^{-\upsilon-1}}{\upsilon-1}+\frac{2ax^{-\upsilon}}{\upsilon-1}\right)\, dx=\\
=-\frac{x^{-\upsilon-1}}{\upsilon-1}e^{-ax^{2}-bx}\quad\left[\upsilon\neq1\right]\end{array}\:\label{eq:generic int Theorem 20-1}\end{equation}

\textbf{Proof }

One may easily take the negative of the indefinite integral of each
line of Theorem 20 with respect to \emph{a} \cite{GR5} (p. 112 Nos.
2.311 and 2.322.1), since there is no error function on the last line,
to complete the proof. $\square$

Theorem 20 with $\upsilon\to3$ gives

\subsection*{Corollary 47}

\begin{equation}
\begin{array}{ll}
\int e^{-ax^{2}-bx}\left(\frac{1}{x^{3}}+\frac{b}{2x^{2}}+\frac{a}{x}\right)\, dx=-\frac{e^{-ax^{2}-bx}}{2x^{2}}\end{array}\:.\label{eq:Theorem 20 n=00003D3}\end{equation}
Theorem 21 with $\upsilon\to3$ gives

\subsection*{Corollary 48}

\begin{equation}
\begin{array}{ll}
\int e^{-ax^{2}-bx}\left(\frac{2}{x^{5}}+\frac{b}{2x^{4}}+\frac{a}{x^{3}}\right)\, dx=-\frac{e^{-ax^{2}-bx}}{2x^{4}}\end{array}\:,\label{eq:a-int of Theorem 20 n=00003D3}\end{equation}
which is not the same as Theorem 20 with $\upsilon\to5$.

Can one step down from Theorem 20 by one power of \emph{x} rather
than two? Indeed,

\textbf{Theorem 22}

\begin{equation}
\begin{array}{ll}
\int e^{-ax^{2}-bx}\left(\frac{\upsilon x^{-\upsilon-1}}{\upsilon-1}+\frac{bx^{-\upsilon}}{\upsilon-1}+\frac{2ax^{1-\upsilon}}{\upsilon-1}\right)\, dx=\\
=-\frac{x^{-\upsilon}}{\upsilon-1}e^{-ax^{2}-bx}\quad\left[\upsilon\neq1\right]\end{array}\:\label{eq:generic int Theorem 22}\end{equation}

\textbf{Proof }

Integrating the resultants with respect to \emph{b} failed for integrals
in Sections 2-5, but for Theorem 20 one may do so since there is no
error function on the last line. One may use \cite{GR5} (p.
112 Nos. 2.311 and 2.322.1) to complete the proof. $\square$

Instead of the integral both lines of Theorem 20 with respect to \emph{b},
one may take the derivative of both lines with respect to \emph{b}, but if we take $\upsilon\to3$,
multiply by $6b^{2}$ and add a constant term multiplied by $b^{2}-6a$,
we see that we have already found this result (\ref{eq:V(0,2,x) Theorem 2}).

Any number of additional theorems may be crafted by taking indefinite
integrals of each line of Theorems 21 and 22 with respect to \emph{a}
or \emph{b}. but we will not expore further down this path.

\section{Appendix C}

Here we provide the Mathematica code for $\Upsilon(p,\, a,\, b,\, x)$
and $V(p,\, n,\, a,\, b,\, x)$. There are no restrictions on the
parameters.

The Mathematica code is fairly straightforward to parse, but I have
also translated the right-hand side of the easier of these, $\Upsilon(p,\, a,\, b,\, x)$,
into Fortran and C as an aid to reprogramming in those languages.
The code for this series produces a table that may be summed using
the Total{[}\%{]} command, which is generally much faster than replacing
Table with the Sum command, since Mathematica spends a great deal of time trying
to find the analytical sum of such series. Note that the code Upsilon{[}p\_,
a\_, b\_, x\_, ifny\_{]} for $\Upsilon(p,\, a,\, b,\, x)$ contains
a fifth parameter {}``ifny'' whose value correponds to the point
at which the infinite series is truncated:

\begin{quote}
Upsilon{[}p\_, a\_, b\_, x\_, ifny\_{]} :=

Table{[}-((a\textasciicircum{}(p/2 - 1/2 - k){*}b\textasciicircum{}(2{*}k){*}
Gamma{[}k - p/2 + 1/2, a{*}x\textasciicircum{}2{]})/(2{*}(2{*}k)!))
+ (a\textasciicircum{}(p/2 - 1 - k){*}b\textasciicircum{}(1 + 2{*}k){*}
Gamma{[}1 - p/2 + k, a{*}x\textasciicircum{}2{]})/(2{*}(1 + 2{*}k)!),
\{k, 0, ifny\}{]}


\vspace{0.3cm}

\noindent Upsilon{[}0,0.11, 0.13, .37, 4{]} $=$ \{-1.72185, -0.0873614,
-0.00173852, -0.0000226837, -2.20537{*}10\textasciicircum{}-7\}.


\end{quote}
\vspace{0.1cm}

Fortran code would loop through and write out results for the expression within the Table command above, translated as

\begin{verbatim} 
        -((a**(p/2. - 0.5 - k)*b**(2*k)*
     -        Gamma(k - p/2. + 0.5,a*x**2))/
     -      (2.*Factorial(2*k))) + 
     -   (a**(p/2. - 1 - k)*b**(1 + 2*k)*
     -      Gamma(1 - p/2. + k,a*x**2))/
     -    (2.*Factorial(1 + 2*k))
\end{verbatim}

\vspace{0.3cm}

\noindent
and the C version of this internal expression is

\begin{verbatim} 
-((Power(a,p/2. - 0.5 - k){*}Power(b,2{*}k){*}
Gamma(k - p/2. + 0.5,a{*}Power(x,2)))/ (2.{*}Factorial(2{*}k))) +
(Power(a,p/2. - 1 - k){*}Power(b,1 + 2{*}k){*} Gamma(1 - p/2. + k,a{*}Power(x,2)))/
(2.{*}Factorial(1 + 2{*}k))) 
\end{verbatim}

Next come the even and odd portions of $\Upsilon(p,\, a,\, b,\, x)$:

\begin{quote}

Upsilone{[}p\_, a\_, b\_, x\_, ifny\_{]} := Table{[}-((a\textasciicircum{}(p/2
- 1/2 - k) b\textasciicircum{}(2 k) Gamma{[}k - p/2 + 1/2, a x\textasciicircum{}2{]})/(
2 (2 k)!)), \{k, 0, ifny\}{]}

\vspace{0.3cm}

Upsilone{[}0,0.11, 0.13, .37, 4{]} $=$ \{-2.30392, -0.102491, -0.00197098,
-0.0000252348, -2.42311{*}10\textasciicircum{}-7\}

\vspace{0.3cm}

Upsilono{[}p\_, a\_, b\_, x\_, ifny\_{]} := Table{[}(a\textasciicircum{}(p/2
- 1 - k) b\textasciicircum{}(1 + 2 k) Gamma{[}1 - p/2 + k, a x\textasciicircum{}2{]})/(
2 (1 + 2 k)!), \{k, 0, ifny\}{]}

\vspace{0.3cm}

Upsilono{[}0,0.11, 0.13, .37, 4{]} $=$ \{0.582077, 0.0151292, 0.000232465,
2.55108{*}10\textasciicircum{}-6, 2.17743{*}10\textasciicircum{}-8\}

\end{quote}

Finally we have the code that one may use to give an analytical or
numerical version of $V(p,\, n,\, a,\, b,\, x)$ for $n=0-5$ and
$p=-n,-n+1,\cdots,\, n-4,\, n-2$. 

\vspace{5.3 cm}

\begin{quote}

V{[}p\_, n\_, a\_, b\_, x\_{]} := Which{[}n == 0, Which{[}p == 0,
( E\textasciicircum{}(b\textasciicircum{}2/(4 a)) Sqrt{[}\textbackslash{}{[}Pi{]}{]}
Erf{[}(b + 2 a x)/(2 Sqrt{[}a{]}){]})/( 2 Sqrt{[}a{]}), p == 1, -(E\textasciicircum{}(-b
x - a x\textasciicircum{}2)/(2 a)) - ( b E\textasciicircum{}(b\textasciicircum{}2/(4
a)) Sqrt{[}\textbackslash{}{[}Pi{]}{]} Erf{[}(b + 2 a x)/(2 Sqrt{[}a{]}){]})/(
4 a\textasciicircum{}(3/2)), p == 2, (E\textasciicircum{}(-b x - a
x\textasciicircum{}2) (b - 2 a x))/( 4 a\textasciicircum{}2) + ((2
a + b\textasciicircum{}2) E\textasciicircum{}(b\textasciicircum{}2/(4
a)) Sqrt{[}\textbackslash{}{[}Pi{]}{]} Erf{[}(b + 2 a x)/(2 Sqrt{[}a{]}){]})/(8
a\textasciicircum{}(5/2)){]}, n == 1, Which{[}p == -1, -(E\textasciicircum{}(-b
x - a x\textasciicircum{}2)/x) - Sqrt{[}a{]} E\textasciicircum{}(b\textasciicircum{}2/(4
a)) Sqrt{[}\textbackslash{}{[}Pi{]}{]} Erf{[}(b + 2 a x)/(2 Sqrt{[}a{]}){]},
p == 1, -((b E\textasciicircum{}(-b x - a x\textasciicircum{}2))/(2
a)) + ((2 a - b\textasciicircum{}2) E\textasciicircum{}(b\textasciicircum{}2/(4
a)) Sqrt{[}\textbackslash{}{[}Pi{]}{]} Erf{[}(b + 2 a x)/(2 Sqrt{[}a{]}){]})/(4
a\textasciicircum{}(3/2)), p == 2, ( E\textasciicircum{}(-b x - a
x\textasciicircum{}2) (b\textasciicircum{}2 - 2 a (1 + b x)))/(4 a\textasciicircum{}2)
+ ( b\textasciicircum{}3 E\textasciicircum{}(b\textasciicircum{}2/(4
a)) Sqrt{[}\textbackslash{}{[}Pi{]}{]} Erf{[}(b + 2 a x)/(2 Sqrt{[}a{]}){]})/(
8 a\textasciicircum{}(5/2)){]}, n == 2, Which{[}p == 0, -((3 E\textasciicircum{}(-b
x - a x\textasciicircum{}2))/x) - ((6 a - b\textasciicircum{}2) E\textasciicircum{}(b\textasciicircum{}2/(4
a)) Sqrt{[}\textbackslash{}{[}Pi{]}{]} Erf{[}(b + 2 a x)/(2 Sqrt{[}a{]}){]})/(2
Sqrt{[}a{]}), p == -2, E\textasciicircum{}(-b x - a x\textasciicircum{}2)
(-(1/x\textasciicircum{}3) - b/x\textasciicircum{}2 + (2 a)/x) + 2
a\textasciicircum{}(3/2) E\textasciicircum{}(b\textasciicircum{}2/(4
a)) Sqrt{[}\textbackslash{}{[}Pi{]}{]} Erf{[}(b + 2 a x)/(2 Sqrt{[}a{]}){]}{]},
n == 3, Which{[}p == 1, 1/(4 a\textasciicircum{}(3/2) x) E\textasciicircum{}(-b
x - a x\textasciicircum{}2) (-2 Sqrt{[} a{]} (30 a + b\textasciicircum{}3
x) - (60 a\textasciicircum{}2 - 12 a b\textasciicircum{}2 + b\textasciicircum{}4)
E\textasciicircum{}((b + 2 a x)\textasciicircum{}2/(4 a)) Sqrt{[}\textbackslash{}{[}Pi{]}{]}
x Erf{[}(b + 2 a x)/(2 Sqrt{[}a{]}){]}), p == -1, 1/x\textasciicircum{}3
E\textasciicircum{}(-b x - a x\textasciicircum{}2) (-5 - 5 b x + 10
a x\textasciicircum{}2 - b\textasciicircum{}2 x\textasciicircum{}2
+ Sqrt{[}a{]} (10 a - b\textasciicircum{}2) E\textasciicircum{}((b
+ 2 a x)\textasciicircum{}2/(4 a)) Sqrt{[}\textbackslash{}{[}Pi{]}{]}
x\textasciicircum{}3 Erf{[}(b + 2 a x)/(2 Sqrt{[}a{]}){]}), p == -3,
E\textasciicircum{}(-b x - a x\textasciicircum{}2) (-(3/x\textasciicircum{}5)
- (3 b)/x\textasciicircum{}4 + (2 a - b\textasciicircum{}2)/x\textasciicircum{}3
+ (2 a b)/x\textasciicircum{}2 - ( 4 a\textasciicircum{}2)/x) - 4
a\textasciicircum{}(5/2) E\textasciicircum{}(b\textasciicircum{}2/(4
a)) Sqrt{[}\textbackslash{}{[}Pi{]}{]} Erf{[}(b + 2 a x)/(2 Sqrt{[}a{]}){]}{]},
n == 4, Which{[}p == 2, E\textasciicircum{}(-b x - a x\textasciicircum{}2)
(-((5 b\textasciicircum{}3)/a) + b\textasciicircum{}5/(4 a\textasciicircum{}2)
- 105/x - (b\textasciicircum{}4 x)/( 2 a)) - ((840 a\textasciicircum{}3
- 180 a\textasciicircum{}2 b\textasciicircum{}2 + 18 a b\textasciicircum{}4
- b\textasciicircum{}6) E\textasciicircum{}(b\textasciicircum{}2/(4
a)) Sqrt{[}\textbackslash{}{[}Pi{]}{]} Erf{[}(b + 2 a x)/(2 Sqrt{[}a{]}){]})/(8
a\textasciicircum{}(5/2)), p == 0, -((5 E\textasciicircum{}(-b x -
a x\textasciicircum{}2) (7 + 7 b x - 14 a x\textasciicircum{}2 + 2
b\textasciicircum{}2 x\textasciicircum{}2))/ x\textasciicircum{}3)
+ ((140 a\textasciicircum{}2 - 20 a b\textasciicircum{}2 + b\textasciicircum{}4)
E\textasciicircum{}(b\textasciicircum{}2/(4 a)) Sqrt{[}\textbackslash{}{[}Pi{]}{]}
Erf{[}(b + 2 a x)/(2 Sqrt{[}a{]}){]})/(2 Sqrt{[}a{]}), p == -2, E\textasciicircum{}(-b
x - a x\textasciicircum{}2) (-(21/x\textasciicircum{}5) - (21 b)/x\textasciicircum{}4
+ (14 a)/x\textasciicircum{}3 - (8 b\textasciicircum{}2)/x\textasciicircum{}3
+ ( 14 a b)/x\textasciicircum{}2 - b\textasciicircum{}3/x\textasciicircum{}2
- (28 a\textasciicircum{}2)/x + (2 a b\textasciicircum{}2)/x) - 2
a\textasciicircum{}(3/2) (14 a - b\textasciicircum{}2) E\textasciicircum{}(b\textasciicircum{}2/(4
a)) Sqrt{[}\textbackslash{}{[}Pi{]}{]} Erf{[}(b + 2 a x)/(2 Sqrt{[}a{]}){]},
p == -4, E\textasciicircum{}(-b x - a x\textasciicircum{}2) (-(15/x\textasciicircum{}7)
- (15 b)/x\textasciicircum{}6 - (6 (-a + b\textasciicircum{}2))/x\textasciicircum{}5
+ ( 6 a b - b\textasciicircum{}3)/x\textasciicircum{}4 + (2 (-2 a\textasciicircum{}2
+ a b\textasciicircum{}2))/x\textasciicircum{}3 - (4 a\textasciicircum{}2
b)/x\textasciicircum{}2 + ( 8 a\textasciicircum{}3)/x) + 8 a\textasciicircum{}(7/2)
E\textasciicircum{}(b\textasciicircum{}2/(4 a)) Sqrt{[}\textbackslash{}{[}Pi{]}{]}
Erf{[}(b + 2 a x)/(2 Sqrt{[}a{]}){]}{]}, n == 5, Which{[}p == 3, E\textasciicircum{}(-b
x - a x\textasciicircum{}2) (-((105 b\textasciicircum{}3)/(2 a)) +
(13 b\textasciicircum{}5)/(4 a\textasciicircum{}2) - b\textasciicircum{}7/(8
a\textasciicircum{}3) - 945/x - (15 b\textasciicircum{}4 x)/(2 a)
+ (b\textasciicircum{}6 x)/(4 a\textasciicircum{}2) - (b\textasciicircum{}5
x\textasciicircum{}2)/( 2 a)) - ((15120 a\textasciicircum{}4 - 3360
a\textasciicircum{}3 b\textasciicircum{}2 + 360 a\textasciicircum{}2
b\textasciicircum{}4 - 24 a b\textasciicircum{}6 + b\textasciicircum{}8)
E\textasciicircum{}(b\textasciicircum{}2/(4 a)) Sqrt{[}\textbackslash{}{[}Pi{]}{]}
Erf{[}(b + 2 a x)/(2 Sqrt{[}a{]}){]})/( 16 a\textasciicircum{}(7/2)),
p == 1, E\textasciicircum{}(-b x - a x\textasciicircum{}2) (-(b\textasciicircum{}5/(2
a)) - 315/x\textasciicircum{}3 - (315 b)/x\textasciicircum{}2 + (630
a)/x - ( 105 b\textasciicircum{}2)/x) + ((2520 a\textasciicircum{}3
- 420 a\textasciicircum{}2 b\textasciicircum{}2 + 30 a b\textasciicircum{}4
- b\textasciicircum{}6) E\textasciicircum{}( b\textasciicircum{}2/(4
a)) Sqrt{[}\textbackslash{}{[}Pi{]}{]} Erf{[}(b + 2 a x)/(2 Sqrt{[}a{]}){]})/(4
a\textasciicircum{}(3/2)), p == -1, E\textasciicircum{}(-b x - a x\textasciicircum{}2)
(-(189/x\textasciicircum{}5) - (189 b)/x\textasciicircum{}4 + (126
a)/x\textasciicircum{}3 - (77 b\textasciicircum{}2)/ x\textasciicircum{}3
+ (126 a b)/x\textasciicircum{}2 - (14 b\textasciicircum{}3)/x\textasciicircum{}2
- (252 a\textasciicircum{}2)/x + (28 a b\textasciicircum{}2)/ x -
b\textasciicircum{}4/x) - Sqrt{[}a{]} (252 a\textasciicircum{}2 -
28 a b\textasciicircum{}2 + b\textasciicircum{}4) E\textasciicircum{}(b\textasciicircum{}2/(4
a)) Sqrt{[}\textbackslash{}{[}Pi{]}{]} Erf{[}(b + 2 a x)/(2 Sqrt{[}a{]}){]},
p == -3, E\textasciicircum{}(-b x - a x\textasciicircum{}2) (-(135/x\textasciicircum{}7)
- (135 b)/x\textasciicircum{}6 + (54 a)/x\textasciicircum{}5 - (57
b\textasciicircum{}2)/x\textasciicircum{}5 + ( 54 a b)/x\textasciicircum{}4
- (12 b\textasciicircum{}3)/x\textasciicircum{}4 - (36 a\textasciicircum{}2)/x\textasciicircum{}3
+ (20 a b\textasciicircum{}2)/x\textasciicircum{}3 - b\textasciicircum{}4/x\textasciicircum{}3
- (36 a\textasciicircum{}2 b)/x\textasciicircum{}2 + (2 a b\textasciicircum{}3)/x\textasciicircum{}2
+ (72 a\textasciicircum{}3)/x - ( 4 a\textasciicircum{}2 b\textasciicircum{}2)/x)
+ 4 a\textasciicircum{}(5/2) (18 a - b\textasciicircum{}2) E\textasciicircum{}(b\textasciicircum{}2/(4
a)) Sqrt{[}\textbackslash{}{[}Pi{]}{]} Erf{[}(b + 2 a x)/(2 Sqrt{[}a{]}){]},
p == -5, E\textasciicircum{}(-b x - a x\textasciicircum{}2) (-(105/x\textasciicircum{}9)
- (105 b)/x\textasciicircum{}8 - (15 (-2 a + 3 b\textasciicircum{}2))/x\textasciicircum{}7
- ( 10 (-3 a b + b\textasciicircum{}3))/x\textasciicircum{}6 + (-12
a\textasciicircum{}2 + 12 a b\textasciicircum{}2 - b\textasciicircum{}4)/x\textasciicircum{}5
+ ( 2 (-6 a\textasciicircum{}2 b + a b\textasciicircum{}3))/x\textasciicircum{}4
- (4 (-2 a\textasciicircum{}3 + a\textasciicircum{}2 b\textasciicircum{}2))/x\textasciicircum{}3
+ ( 8 a\textasciicircum{}3 b)/x\textasciicircum{}2 - (16 a\textasciicircum{}4)/x)
- 16 a\textasciicircum{}(9/2) E\textasciicircum{}(b\textasciicircum{}2/(4
a)) Sqrt{[}\textbackslash{}{[}Pi{]}{]} Erf{[}(b + 2 a x)/(2 Sqrt{[}a{]}){]}{]}{]}

\end{quote}

\vspace{0.3cm}

\begin{quote}

V{[}0, 0,0.11, 0.13, .37{]} $=$ 0.965734

\vspace{0.3cm}
 V{[}1, 0,0.11, 0.13, .37{]} $=$ -4.83791

\vspace{0.3cm}
 V{[}2, 0,0.11, 0.13, .37{]} $=$ 5.66958

\vspace{0.3cm}
 V{[}-1, 1,0.11, 0.13, .37{]} $=$ -2.74974

\vspace{0.3cm}
 V{[}1, 1,0.11, 0.13, .37{]} $=$ 0.336806

\vspace{0.3cm}
 V{[}0, 2,0.11, 0.13, .37{]} $=$ -8.23291
 
  \vspace{0.3cm}
 V{[}-2, 2,0.11, 0.13, .37{]} $=$ -18.8204

\vspace{0.3cm}
 V{[}1, 3,0.11, 0.13, .37{]} $=$ -41.1588

\vspace{0.3cm}
 V{[}-1, 3,0.11, 0.13, .37{]} $=$ -94.1483

\vspace{0.3cm}
 V{[}-3, 3,0.11, 0.13, .37{]} $=$ -421.855

\vspace{0.3cm}
 V{[}2, 4,0.11, 0.13, .37{]} $=$ -288.093

\vspace{0.3cm}
 V{[}0, 4,0.11, 0.13, .37{]} $=$ -659.177

\vspace{0.3cm}
 V{[}-2, 4,0.11, 0.13, .37{]} $=$ -2953.3

\vspace{0.3cm}
 V{[}-4, 4,0.11, 0.13, .37{]} $=$ -15468.2

\vspace{0.3cm}
 V{[}3, 5,0.11, 0.13, .37{]} $=$ -2592.75

\vspace{0.3cm}
 V{[}1, 5,0.11, 0.13, .37{]} $=$ -5933.29

\vspace{0.3cm}
 V{[}-1, 5,0.11, 0.13, .37{]} $=$ -26581.3

\vspace{0.3cm}
 V{[}-3, 5,0.11, 0.13, .37{]} $=$ -139221.

\vspace{0.3cm}
 V{[}-5, 5,0.11, 0.13, .37{]} $=$ -792320.

\end{quote}


\end{document}